\documentclass[12pt]{article}
\usepackage{amssymb}
\usepackage{enumerate}
\usepackage{amsfonts}
\usepackage{amsmath}
\usepackage{amsthm}
\usepackage{psfrag} 
\usepackage{epsfig}
\usepackage{subfigure}
\usepackage{graphicx}
\usepackage{color}
\usepackage{amssymb}
\usepackage{epstopdf}
\usepackage{amsbsy}
\usepackage{latexsym}
\usepackage{moreverb}                      
\usepackage{textcomp}

\usepackage{multicol}
\usepackage{algorithm}
\usepackage{algorithmic}

\newcommand{\lbeq}[1]{{\label{OR:eq:#1}}}
\newcommand{\be}[1]{\begin{equation} \lbeq{#1}}
\newcommand{\ee}{\end{equation}}
\newcommand{\beno}{\begin{equation*}}
\newcommand{\eeno}{\end{equation*}}

\newtheorem{remark}{Remark}
\newtheorem{definition}{Definition}
\newtheorem{theorem}{Theorem}

\newcommand{\ubar}[1]{\text{\b{$#1$}}}

\usepackage[top=2.6cm,bottom=2.6cm,left=2cm,right=2cm]{geometry}

\usepackage[T1]{fontenc}
\usepackage[utf8]{inputenc}

\usepackage[affil-it]{authblk}
\usepackage{lmodern}


\usepackage{sectsty}
\sectionfont{\fontsize{13}{15}\selectfont}
\subsectionfont{\fontsize{12}{15}\selectfont}
\begin{document}

\title{A Fuzzy-Stochastic Multiscale Model for Fiber Composites\\
{\large{A one-dimensional study}}}

\author[1]{Ivo Babu{\v s}ka\thanks{babuska@ices.utexas.edu}}
\author[2]{Mohammad Motamed\thanks{motamed@math.unm.edu}}

\affil[1]{Institute for Computational Engineering and Sciences, The University of Texas at Austin, USA}
\affil[2]{Department of Mathematics and Statistics, The University of New Mexico, Albuquerque, USA}

\maketitle

\begin{abstract}
We study mathematical and computational models for computing the deformation of fiber-reinforced cross-plied laminates due to external forces. This requires an understanding of both micro-structural effects and different sources of uncertainty in the problem. We first show that the uncertainties in the problem are of both statistical (aleatoric) and systematic (epistemic) types and that current multiscale stochastic models, such as stationary random fields, which are based on {\it precise probability} theory, are not capable of correctly characterizing uncertainty in fiber composites. Next, we motivate the applicability of models based on {\it imprecise uncertainty} theory and present a novel fuzzy-stochastic model, which can more accurately describe uncertainties in fiber composites. The new model is constructed by combining stochastic fields and fuzzy variables through a simple  calibration-validation approach. Finally, we construct a global-local multiscale algorithm for efficiently computing output quantities of interest. The method aims at approximating required quantities, such as displacements and stresses, in regions of relatively small size, e.g. hot spots or zones. The algorithm uses the concept of representative volume elements and computes a global solution to construct a local approximation that captures the microscale features of the solution. The results are based on and backed by real experimental data. 
\end{abstract}


\section{Introduction}

\begin{figure}[!h]
\begin{multicols}{2}
\centering
\includegraphics[width=5.5cm,height=4cm]{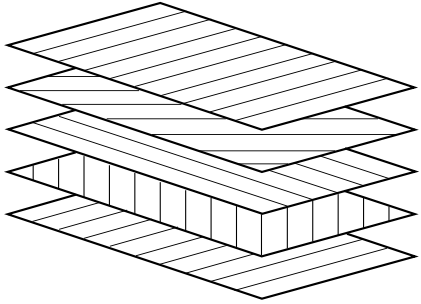}\\
\vskip .03cm
(a) Plies with fibers in different directions

\includegraphics[width=.9\linewidth]{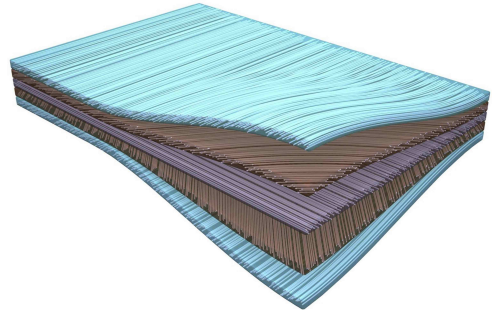}\\
\vskip .15cm
(b) Plies are stacked to form a laminate 
\end{multicols}
\vskip -0.3cm
\caption{A cross-plied fiber composite laminate made by stacking five plies. Different plies contain many unidirectional fibers in different directions aligned in a matrix.}
\label{laminate}
\end{figure}

Fiber-reinforced composite materials are widely used in aerospace, marine, and automotive industries. They consist of stiff fibers in a matrix which is less stiff. In composites with unidirectional fibers, a large number of long unidirectional fibers are aligned in a thin ply. To achieve high stiffness, a few plies are stacked together, each having fibers oriented in a certain direction. Such a stack is termed a cross-plied laminate; see Figure \ref{laminate}. 

Fiber composite materials are designed and fabricated based on their response to external forces. Two principal responses are deformation and fracture. We therefore need a through understanding of stress distributions and damage and failure mechanisms, such as the initiation and propagation of matrix cracks, fiber breakage, and fiber-matrix interface splitting. 
In a non-testing environment, a through understanding of the mechanical behavior of fiber composites must be obtained from accurate and viable mathematical and computational models. 
The overall mechanical behavior is very complex and requires an understanding of both micro-structural effects and variability/uncertainty in the manufacturing process. A mathematical model must therefore consider and include two components: {\it multiple scales} and {\it uncertainty}. On the one hand, the problem involves multiple length scales ranging from the diameter of fibers (on the order of $10^{-6}\,$m) to the laminate thickness and length (on the order of $10^{-2}-10^{1}\,$m). On the other hand, the model is subject to uncertainty, due to the random character of the size, location, and distribution of fibers and the intrinsic variability in materials properties and fracture parameters. High levels of confidence in the predictions require an understanding of the uncertainties in the multiscale model. This understanding can be obtained by a process called uncertainty quantification (UQ) taking into account the systematic coupling and interaction of the two involved components, i.e. multiple scales and uncertainty.

A chronological development of mathematical and computational tools for solving multiscale problems (such as fiber composites) may in general be classified into four groups:
\begin{itemize}
\item[1.] Mathematical theory of homogenization: The term ``{\it homogenization}'' was first coined by Babu{\v s}ka \cite{Babuska1:74, Babuska2:75}. Homogenization is an analytical approach to replace the multiscale problem with heterogeneous coefficients by an equivalent problem with homogeneous coefficients, known as the homogenized problem. Homogenization theory is well studied in the case of periodic or locally periodic microscale coefficients \cite{Bensoussan_etal:78}. For the more general case of non-periodic coefficients, homogenization theory is carried out by studying the G-convergence \cite{Spagnolo:75} and H-convergence \cite{Murat_Tartar:78} of solution operates. 

\item[2.] Multiscale numerical methods: Numerical methods that approximate the highly oscillatory solution of a multiscale problem by solving an effective problem and including local microscale oscillations instead of  directly solving the original problem are called multiscale numerical methods or numerical homogenization methods. They were first introduced and studied by Babu{\v s}ka \cite{Babuska3:75,Babuska_I:76,Babuska_II:76,Babuska_III:77}, 
where the finite element method was employed to solve the homogenized problem and some microscale features were added as correctors to the homogenized solution. Since then, many multiscale methods have been proposed, including generalized finite element methods \cite{Babuska_Osborn:83,Matache_etal:2000,Babuska_Lipton:11b}, 
variational multiscale methods \cite{Hughes:95,Hughes:98}, multiscale finite element methods \cite{Efendiev_Hou:09}, 
projection-based numerical homogenization \cite{Dorobantu_Engquist:98,Engquist_Runborg:02b}, 
heterogeneous multiscale methods \cite{E_Engquist:03,E_etal:07}, 
and equation-free methods \cite{Kevrekidis_etal:03}. We also refer to \cite{Hou:03,E:11,Engquist_Souganidis:08} for detailed discussions on a wide range of multiscale methods. All these techniques usually seek the approximate solution everywhere inside the computational domain. Another class of multiscale methods is based on a global-local approach, where the solution or other quantities of interest (QoIs) are required inside a relatively small subdomain. In these methods, a global (homogenized) solution is employed to recover the microscale solution inside the local subdomain; see e.g. \cite{Oden1,Oden2,Babuska_Lipton:11a}. 

\item[3.] Multiscale methods in engineering: There is a vast literature on multiscale methods in engineering, especially in the field of materials science. These methods have been developed based on the same ideas and principles as those developed in the applied mathematics community and mentioned above. They have been proposed to account for microstructural heterogeneity, for instance in complex materials such as composites and porous structures. They include unit cell methods or representative volume element (RVE) approach, multi-level approaches using the finite element method and the Voronoi cell finite element method, and continuous-discontinuous homogenization; see for instance \cite{Kanoute_etal:09,Geers_etal:10,Ladeveze_Lubineau_Marsal:06,Soutis_Beaumont:05,Ghosh:11,Fish_etal:99,Feyel_Chaboche:03}. These multiscale methods treat both linear and non-linear mechanical and thermomechanical responses of complex materials.

\item[4.] Probabilistic treatment of multiscale problems: The literature on UQ for multiscale problems is rather sparse and focuses more on stochastic models. Stochastic homogenization \cite{Souganidis:99,Bourgeata_Piatnitski:04,Blanc_LeBris:10,Gloria_Otto:11} 
can be considered as a generalization of classical homogenization. Theoretical aspects of stochastic homogenization is well studied in the case of stationary and ergodic random fields. However, numerical approaches based on stochastic homogenization are cumbersome and not well studied, particularly because the homogenized problem is set on the whole space, not on a finite cell. Analogous to the case of deterministic problems, a variety of stochastic multiscale methods have been proposed within the framework of variational multiscale methods, multiscale finite element methods, the heterogeneous multiscale methods, and the global-local approach; see e.g. \cite{Asokan_Zabaras:06,Ganapathysubramanian_Zabaras:07,Aarnes_Efendiev:08,Ma_Zabaras:11,LeBris_etal:14,BMT:14}. Also in engineering community, an increasing number of papers are attempting to address UQ in multiscale simulations; see e.g. \cite{XU_Graham-Brady:05,Guilleminota_etal:09,Liu_etal:10,Liu_etal:13,Papadrakakis_Stefanou:14}. 

\end{itemize}


In the particular case of composite materials, a majority of multiscale models are deterministic and hence do not include and treat the present uncertainties; see e. g. \cite{Ladevez_etal:03,Ladaveze_etal:06,Ladaveze_etal:12} in addition to references in item 3 above. 
Such models are not capable of fully describing the mechanical behavior of composites, partially due to the ignorance of uncertainty which is an important component that must be included in the model. 
More recently, there have been efforts to include and characterize uncertainty and to propose stochastic models for fiber composites, see e. g. \cite{Ibnabdeljalil_Phoenix:95,Babuska_etal:1999} in addition to references in item 4 above. 
Such models describe and treat uncertainty by {\it precise probabilities} \cite{Grigoriu}, where the model input parameters, i.e. materials properties such as the modulus of elasticity, are described by (often stationary Gaussian or log-normal) random fields. Despite recent advances with stochastic multiscale models and UQ methodologies, five decades after the pioneering work of Kachanov \cite{Kachanov:58,Kachanov:86} on continuum damage mechanics, one question still remains open: are there pertinent models and tools for the sound investigation of composite responses to micro-defects? As Rohwer \cite{Rohwer:14} has recently remarked, ``A fully satisfying model for describing damage and failure of fiber composites is not yet available. Consequently, for the time being, a real test remains the authentic way to secure structural strength.''

A main objective of the present work is to show the deficiencies of stochastic models, based on precise probability, in accurately predicting fiber composite responses. We are particularly concerned with the deformation of fiber composites due to external forces. Since realistic mathematical and computational models must be designed based on and backed by real experimental data, we consider a small piece of a real fiber composite plate, taken from \cite{Babuska_etal:1999}, consisting of four plies and containing $13688$ unidirectional fibers with a volume fraction of $63 \%$. The real measured data include the size, location, and distribution of fibers, which are obtained by an optical microscope. Using the available real measurements, as a prototype of fiber distributions in fiber composites, we show that the current  stochastic models, such as stationary random fields, are not capable of correctly characterizing uncertainty in fiber composites. Instead, we study and show the applicability of {\it imprecise uncertainty} models for fiber composites. We propose a new model constructed by combining stochastic fields and fuzzy variables \cite{Zadeh:65} through a calibration-validation approach. We finally present a numerical method for propagating uncertainty through the proposed model and predicting output QoIs. The numerical method uses the concept of RVEs and is based on a global-local approach, where a global solution is used to construct a local solution that captures the micro scale features of the solution. For simplification and to motivate and establish the main concepts of the proposed model, we consider a one-dimensional problem. Our one-dimensional studies show the deficiencies of precise probability and motivate the importance and applicability of imprecise uncertainty for modeling the responses of materials with a microstructure. The present work is a preparation for studying fiber composites in two and three dimensions, which will be presented elsewhere.

The main contributions of this paper include: (1) showing the deficiency of precise probability models to the reliable prediction of fiber composite responses; (2) motivating the applicability of imprecise uncertainty models and constructing a novel fuzzy-stochastic model for fiber composites; and (3) developing a global-local numerical method for an efficient propagation of uncertainty through the multiscale model and computing the output QoIs.

The rest of the paper is organized as follows. In Section \ref{sec:PS} we present the available real data and the mathematical formulation of the problem. We then briefly address different uncertainty models for characterizing uncertainty. In Section \ref{sec:statistic} we perform statistical analysis of the data, including their statistical moments and correlation. We present the basic concepts of fuzzy set theory in Section \ref{sec:fuzzy}. In Section \ref{sec:modeling}, we present and discuss the construction of the fuzzy-stochastic model. We propose a global-local algorithm in Section \ref{sec:global-local}. Finally, we summarize our conclusions and outline future works in Section \ref{sec:CON}.

\section{Problem Statement}
\label{sec:PS}

Reliable mathematical and computational models for predicting the response of fiber composites due to external forces must be designed based on and backed by {\it real experimental data}. In this section, we first present the real data that is used throughout this work. We then formulate a simple one-dimensional problem describing the deformation of fiber composites. Finally, we briefly address different models for characterizing uncertainty in the problem. 


\subsection{Real data}
\label{sec:real_data}

\begin{figure}[!h]
\centering
\subfigure{\includegraphics[width=9cm,height=4cm]{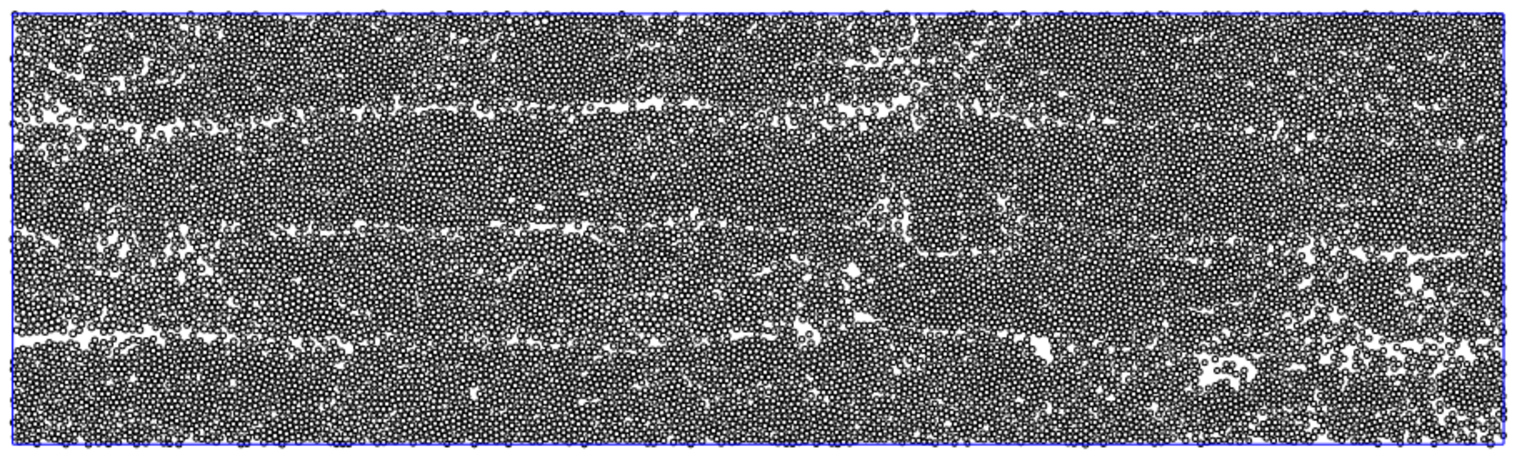}}
       \hskip 1cm
       \subfigure{\includegraphics[width=3.4cm,height=3.6cm]{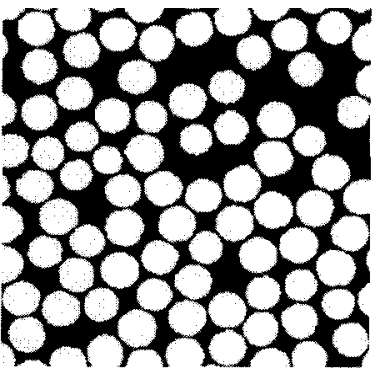}}
\vskip -.3cm
\caption{{\it Left}: A $1.7 \times 0.5 \, \text{mm}^2$ rectangular orthogonal cross section of a small piece of a fiber composite laminate consisting of four uni-directional plies containing 13688 fibers with a volume fraction of $63 \%$. {\it Right}: A binary image of a small part of the whole micrograph.}
\label{composite_fig}
\end{figure}

The real data that we use are obtained from a small piece of a HTA/6376 carbon fiber-reinforced epoxy composite plate \cite{Babuska_etal:1999,BMT:14} with a rectangular cross section of size $1.7 \times 0.5 \, \text{mm}^2$, 
and consisting of four plies containing $13688$ unidirectional fibers with a volume fraction of $63 \%$. Fiber diameters vary between $4 \,\mu$m to $10 \,\mu$m. Figure \ref{composite_fig} shows a map of the size and position of fibers in an 
orthogonal cross section of the composite obtained by an optical microscope. 
In the present work, this particular map serves as a prototype of fiber distributions in fiber composites. 

The Young's modulus of elasticity and Poisson's ratio of the fiber composite under consideration are given in Table \ref{table_material}.
\begin{table}[!htb]
\centering 
\caption{Material constants for the composite under consideration.}
\vskip .2cm
{
\renewcommand{\arraystretch}{1.3}
\begin{tabular}{ccc}
\hline
composite phases &  $\qquad$ $a$  $\qquad$ & $\qquad$ $\nu$ $\qquad$\\
\hline 
fiber &  24 [GPa] & 0.24 \\
matrix &  3.6 [GPa]  & 0.3 \\
\hline
\end{tabular}
}
\label{table_material}
\end{table}


\subsection{Mathematical formulation: a one-dimensional problem}
\label{sec:PS-math}

The deformation of elastic materials is given by the elastic partial differential equations (PDEs) in three dimensions. In the particular case of plane strain, where the length of structures in one direction is very large compared to the size of structures in the other two directions, the problem may approximately be reduced to a two-dimensional problem. In the present work, however, for simplification and to motivate and establish the main concepts of the proposed model, we consider a one-dimensional problem. Problems in higher dimensions will be presented elsewhere. 

We consider the elastic equation with homogeneous Dirichlet and non-homogeneous Neumann boundary conditions in one dimension:
\begin{subequations}\label{BVP}
\begin{gather}
\frac{d}{dx} \bigl( a(x) \, \frac{du}{dx}(x) \bigr) = f(x), \hskip 1cm  x \in [0,1], \label{PDE}\\
u(0)= 0, \quad a(1) \,\frac{du}{dx} (1)= 1, \label{BCs}
\end{gather}
\end{subequations}
where $x$ is location, $u(x)$ is the displacement, $a(x)$ is the modulus of elasticity of the composite, and $f(x)$ is a force term, given for instance by,
\begin{eqnarray*}
f(x) = \left\{ \begin{array}{l l}
2 \, & \, \, x \in [0 , 0.5),\\
0 \, & \, \, x \in [0.5 , 1].
\end{array} \right.
\end{eqnarray*}
Here, the unit of length is assumed to be meter $[\text{m}]$, and the unit of force (both external force $f$ and the boundary force) and modulus of elasticity $a$ is assumed to be giga Pascal $[\text{GPa}]$. 

The main goal of computations is to obtain composite deformations due to external forces. We therefore need to solve problem \eqref{BVP} and compute displacements $u(x)$, stresses $a(x) \, u'(x)$, and/or other QoIs (for example functionals of the displacement $u(x)$). We note that in the one-dimensional model problem considered here, stresses are smooth functions and do not oscillate with the small scale of fiber sizes. We therefore set our main goal as the prediction of the solution at a given point, say at $x_0=0.75$. We then introduce the QoI:
\begin{equation}\label{QOI}
{\mathcal Q} = u(x_0), \qquad x_0 = 0.75.
\end{equation}
Here, the solution to the problem \eqref{BVP} is analytically given by
\begin{equation}\label{exact_sol}
u(x) = \int_{0}^{x} a^{-1}(\xi) \, F(\xi) \, d\xi, \qquad F(\xi) = \int_{0}^{\xi} f(\tau) \, d\tau = \left\{ \begin{array}{l l}
2 \, \xi \, & \, \, \xi \in [0 , 0.5),\\
1 \, & \, \, \xi \in [0.5 , 1].
\end{array} \right.
\end{equation}

To obtain the solution and the QoI in \eqref{QOI}, we need to know the parameter $a(x)$, which describes the mechanical property of the composite. Therefore, we first need to characterize the modulus of elasticity $a(x)$, which is directly given by the size and position of fibers in the matrix and by the modulus of elasticity of fibers and matrix, listed in Table \ref{table_material}.


The parameter $a(x)$ needs to be modeled based on {\it real data} provided by the manufacturer. Moreover, as mentioned in introduction, two major difficulties arise in modeling the material parameter. First, the parameter is highly oscillatory, because there are a large number of small fibers distributed in the matrix. A direct representation is therefore not practical due to computational power limitations and unavailability of data for the whole structure. Secondly, uncertainty---due to various sources---must be included in the expression of the parameter. 


\subsection{On uncertainty characterization and stochastic models}
\label{sec:PS-uncertainty}

Current models for characterizing the parameter $a(x)$ are based on stochastic representation of uncertainty. 
%
%
In practice, the parameter is assumed to be a stationary (and often Gaussian or lognormal) random field with known marginal probability distributions. 
However, as we will show later, stochastic models in general, and stationary fields in particular, are not able to correctly represent the uncertainty in materials properties. Other models need to be considered.
%
%
%
%
%
%
In general, the uncertainty in $a$ can be classified into two categories:
\begin{itemize}

\item[{\bf I}.] {\it Precise or deterministic uncertainty}. This uncertainty is associated with an identifiable but unknown or uncertain source. In this case, we exactly know the uncertainty in the parameter. Two common models include: (1) stochastic models, where $a$ is modeled by random variables or random fields with known and crisp distributions and covariances; and (2) fuzzy models, where $a$ is modeled by fuzzy variables or fuzzy fields with known and crisp membership functions; see Section \ref{sec:fuzzy} for the definition of membership function.

\item[{\bf II}.] {\it Imprecise or non-deterministic uncertainty}. In this case, we do not exactly know the uncertainties in the parameter. The uncertainty characterization can be carried out by different combinations of the above two models. A few examples include: (1) mixed stochastic models, where $a$ is modeled by a random field, whose mean and variance are random variables, (2) mixed fuzzy models, where $a$ is modeled by a fuzzy field, whose mean is a fuzzy number; and (3) fuzzy-stochastic models, where $a$ is modeled by a random field, whose mean and variance are fuzzy numbers.
\end{itemize}

Precise uncertainty, and particularly stochastic models, can be used to characterize uncertainty in physical systems in the presence of abundant, accurate data. In many real applications such as fiber composites, however, the available data are scarce and imprecise. In such cases, the parameters cannot be accurately characterized by precise uncertainty models. Moreover, as we will show in details in the next section, the materials properties are not stationary fields. The non-stationarity is also visible from the map of fibers in Figure \ref{composite_fig}. 
%
%
In the rest of the paper, we first show that stochastic models may not be capable of correctly characterizing uncertainty in fiber composites. We then propose a combined fuzzy-stochastic model for the uncertain microstructure parameter $a$.

\section{Statistical analysis}
\label{sec:statistic}
In this section we perform statistical analysis on the data. We first describe data collection procedures. 
We then study statistical moments and correlation length of the parameter field.

\subsection{Data collection}
\label{sec:data_collection}

\noindent
{\bf A binary data map.} 
The available data contains the center's location and area of all 13688 fibers in the orthogonal cross section  $D \subset {\mathbb R}^2$ of the composite shown in Figure \ref{composite_fig}. We first discretize the two-dimensional domain $D$ into a uniform mesh consisting of $500 \times 1700$ square elements (pixels) of size $1 \times 1 \, \mu {\text{m}}^2$. We then construct a binary data structure, where the presence or absence of fiber at every pixel is marked by 1 or 0, respectively, assuming that fibers are perfectly circular. Using the binary structure generated above, the material elasticity moduli are recovered. For instance, to compute Young's modulus, at each pixel with labels 1 or 0, we put $a = a_{\text{fiber}} = 24$ [Gpa] or $a = a_{\text{matrix}} = 3.6$ [Gpa], respectively. 
We note that significant experiments have been performed for validating the binary mapping procedure and the size of square elements or pixels. It was found that a size of $1 \, \mu \text{m}$ provides the desired accuracy without being too computationally expensive. In particular, the binary map with elements of size 1 micron is accurate within $1 \%$ in predicting the overall volume fraction obtained by an analytic (and expensive) approach. 

%
%

\begin{figure}[!h]
\psfrag{amh10}[c][][0.8]{$\text{1D sample \,} a_m(x_j)$}
\psfrag{binary_micron_data}[c][][0.8]{binary-micron data}
\psfrag{harmonic_average}[c][][0.8]{harmonic average}
\psfrag{j1}[c][][0.8]{$j=1$}
\psfrag{j2}[c][][0.8]{$j=2$}
\psfrag{jend}[c][][0.8]{$j=170$}
\vskip .3cm
\centering
\includegraphics[width=10cm,height=5.3cm]{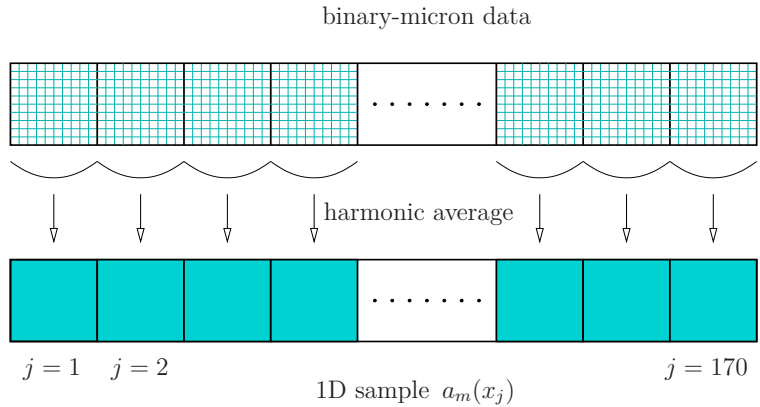}
\vskip -.3cm
\caption{A schematic representation of deriving 1D data samples (bottom) from the binary data given on the $1 \times 1 \, \mu$m pixels (top). On each element $j$, we take the harmonic average of $a$ on its $10 \times 10$ pixels.}
\label{data_collect}
\end{figure}

\medskip
\noindent
{\bf One-dimensional data collection.} 
To obtain one-dimensional data samples from the binary data map, we proceed as follows. First, we divide the rectangular cross section of the composite $D$ into $M=50$ thin horizontal strips  (or bars) of width $h=10 \, \mu$m and length $1700 \, \mu$m along the $x$-axis. This gives us $M$ thin bars labeled $m=1,\dotsc,M$. By the binary-micron map constructed above, the modulus of elasticity in each bar is given on a mesh of $10 \times 1700$ one-micron pixels. Next, each bar $m$ is divided into 170 square elements of size $10 \, \mu$m, labeled $j=1, \dotsc, 170$; see Figure \ref{data_collect} (top). Each element $j$ in a bar contains $10 \times 10$ pixels. On each element $j$, we take the harmonic average over its $10 \times 10$ pixels and compute a one-dimensional modulus of elasticity $a_m(x_j)$. We repeat the process for all 170 elements along $x$ axis and compute $\{ a_m(x_j) \}_{j=1}^{170}$. This is considered as a sample of the uncertain parameter $a(x)$; see Figure \ref{data_collect} (bottom). Eventually we compute and collect all $M=50$ samples by repeating the same steps for all thin bars and obtain $M$ one-dimensional  samples $\{ a_{m}(x_j) \}_{m=1}^M$ at the discrete points $\{ x_j \}_{j=1}^{170}$. 

We note that the parameter $a(x)$ involves multiple scales and is highly oscillatory. The small variations in $a(x)$ are of the size of fibers ($ \sim10 \, \mu$m) and are much smaller ($10^5$ times) than the length of the whole domain ($\sim 1$ m). We obtain the samples on elements of size $10 \, \mu$m in order to fully resolve the micro scales. The harmonic averaging is also motivated by one-dimensional periodic homogenization. We emphasize however that since the size of elements are approximately of the size of fibers, the samples obtained here are not homogenized, and they correspond to a full resolution of the micro scale. We will use homogenization in Section 6, where a global-local approach will be presented to solve the multiscale problem.

\medskip
\noindent
{\bf Bootstrapping.} 
To this end, we have $M=50$ discrete and highly oscillatory samples of $a(x)$ given on $x \in [0, 1.7]$ mm. Each sample corresponds to a one-dimensional bar of length $L=1.7$ mm. 
In practice, we may need more samples ($M>50$) of bars with longer lengths ($L > 1.7$ mm), particularly because the length of the composite plate (1 m) is much larger than the length of the available data (1.7 mm). This can be done by bootstrapping technique \cite{Bootstrap1,Bootstrap3}. 
For example, suppose that from the original data set (i.e. 50 samples of length 1.7 mm), we want to generate $M=100$ bootstrap samples of length $L=10$ mm. A bootstrap sample is generated by randomly sampling $n=\lceil 10/1.7 \rceil = 6$ times, with replacement, from the original set and putting them next to each other. This generates a bar of length $6*1.7 = 10.2$ mm. The extra $0.2$ mm is then removed to obtain a sample of length 10 mm. We repeat this process and obtain the rest of $M-1$ bootstrap samples of length 10 mm. This procedure can be performed for any given number of samples $M$ and any desired length $L$. The new $M$ bootstrap samples of length $L$, denoted by $\{ a_m(x_j) \}_{m=1}^M$ at $N_x = 1+ L [\mu \text{m}] / 10$ discrete points $\{ x_j \}_{j=1}^{N_x} \in [0, L]$ will be used as the data for statistical analysis.

\subsection{First four statistical moments}
\label{sec:moments}

\begin{figure}[!h]
\vskip -.3cm
  \begin{center}
        \subfigure[sample mean]{\includegraphics[width=8.2cm,height=6.5cm]{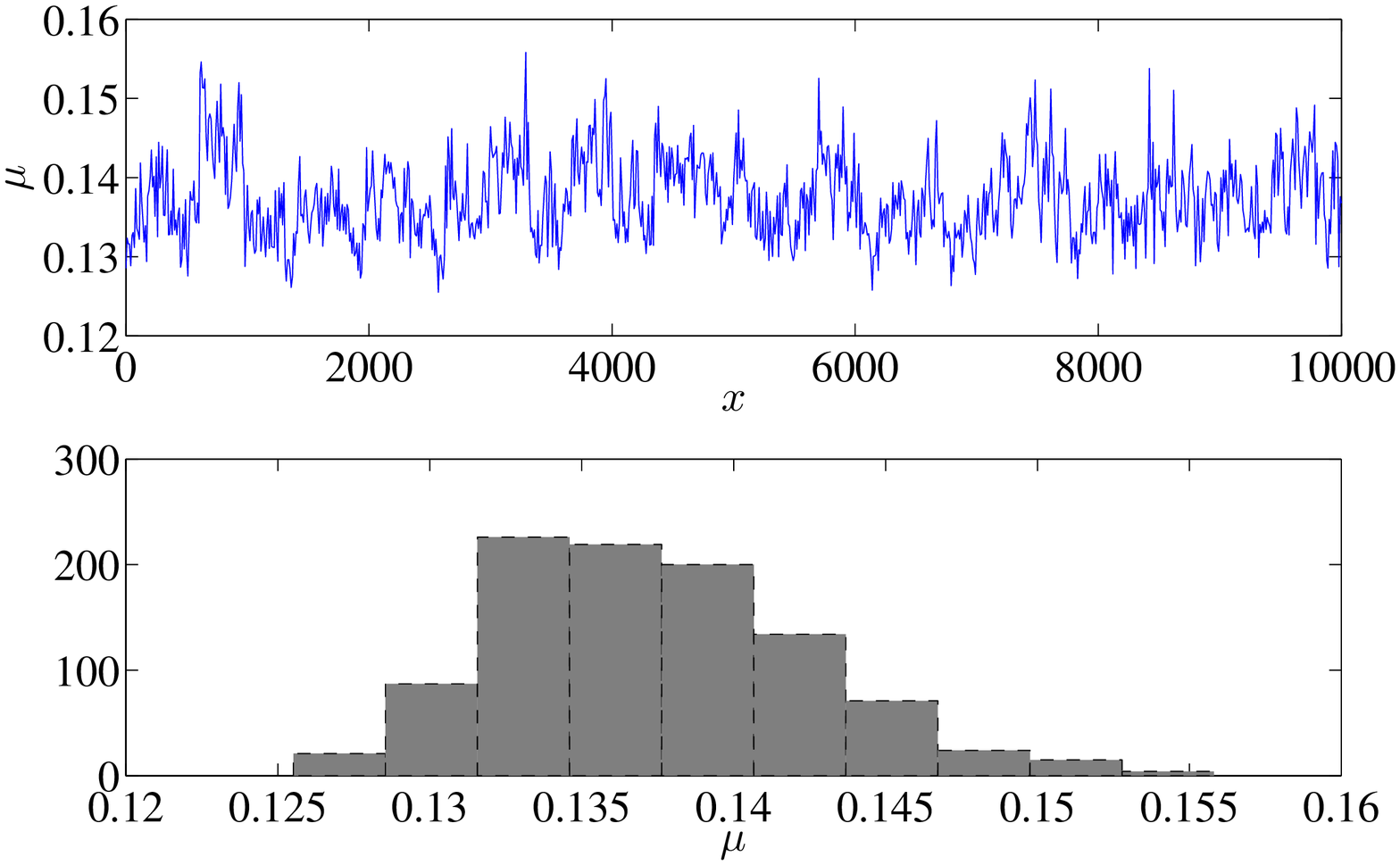}}
        \subfigure[sample standard deviation]{\includegraphics[width=8.2cm,height=6.5cm]{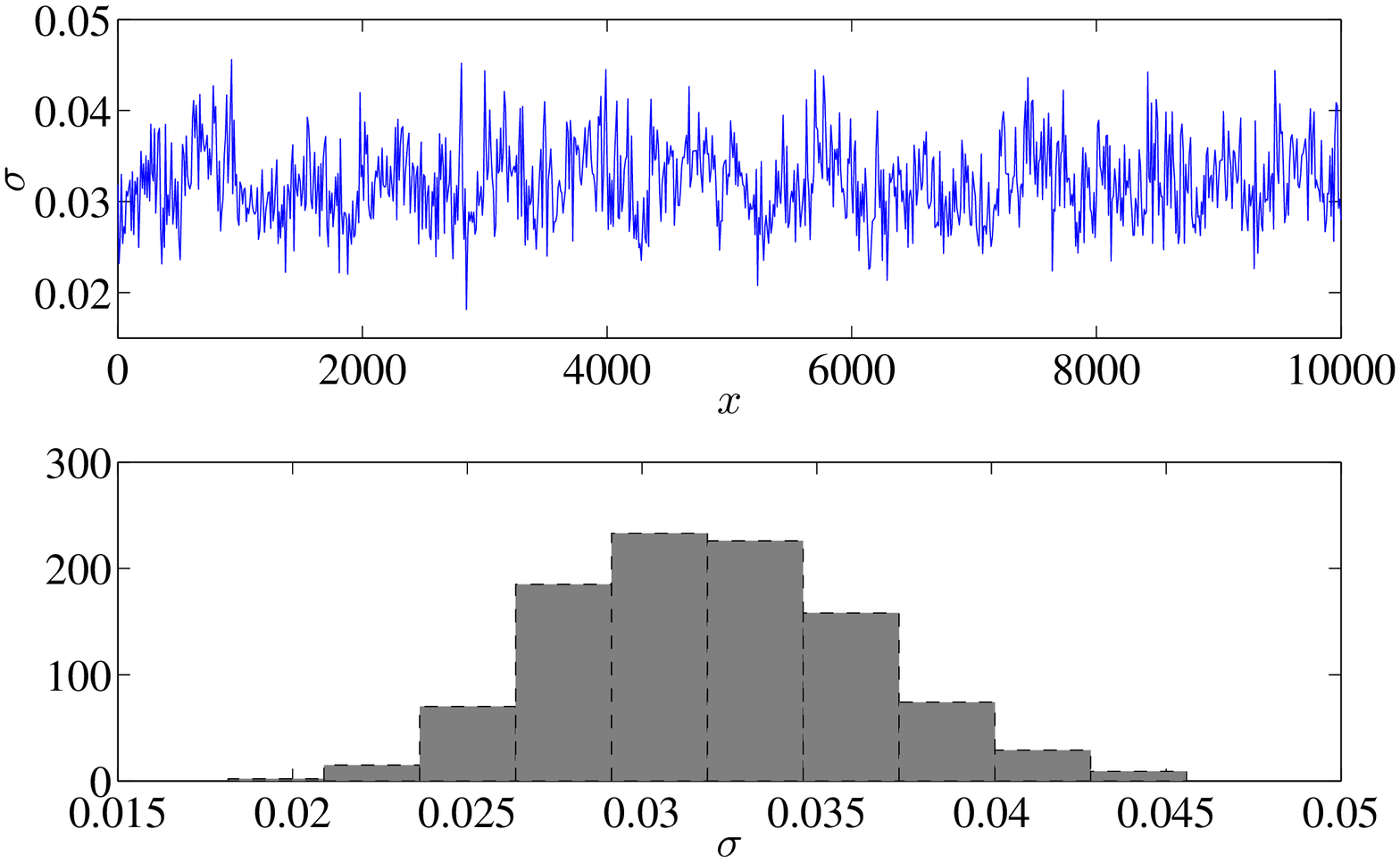}}\\
        \subfigure[sample skewness]{\includegraphics[width=8.2cm,height=6.5cm]{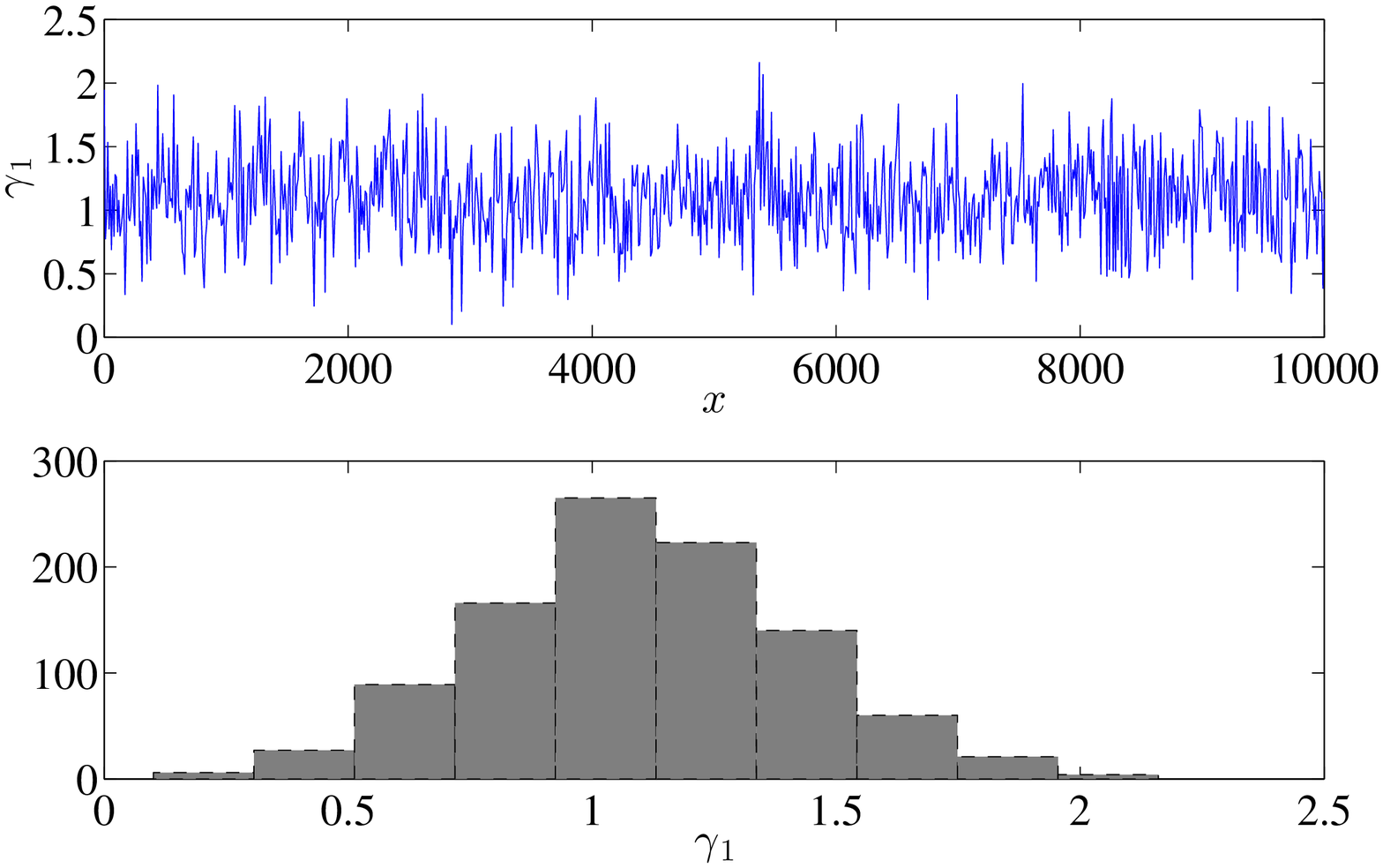}}
        \subfigure[sample excess kurtosis]{\includegraphics[width=8.2cm,height=6.5cm]{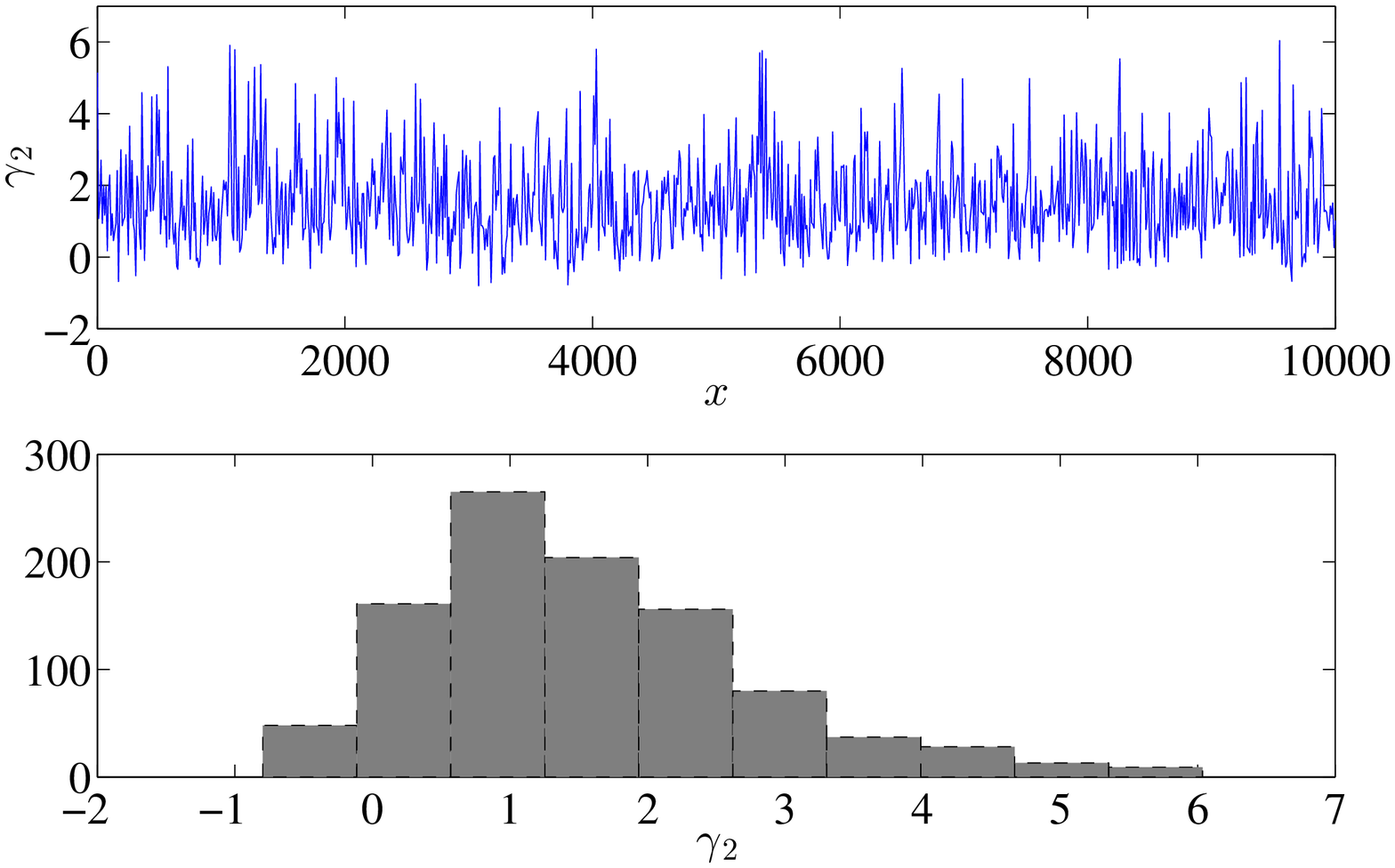}}
\vskip -0.3cm        
\caption{Sample moments of the field $b(x) = a^{-1}(x)$ versus $x$ and their histograms.}
\label{mean_std}
  \end{center}
  \vskip -.5cm
\end{figure}

Motivated by the analytical form of the solution \eqref{exact_sol}, we will perform statistical analysis on the reciprocal of $a$ instead of $a$ itself. We therefore set $b(x) = a^{-1}(x)$ and carry out the analysis using $M$ discrete bootstrap samples $\{ b_{m}(x_j) \}_{m=1}^M$ of length $L$, obtained in Section \ref{sec:data_collection}, where $b_{m}(x_j) = a_m^{-1}(x_j)$. As an example in this section, we choose $M=100$ and $L=10$ mm, and approximate the first four moments of the field $b(x)$ at each discrete point $\{ x_j \}_{j=1}^{N_x}$ using $M=100$ samples $\{ b_{m}(x_j) \}_{m=1}^M$ of length $L=10$ mm:
%
%
$$
\mu(x_j) = \frac{1}{M} \sum_{m=1}^{M} b_{m}(x_j), \qquad  \sigma^2(x_j) =\frac{1}{M} \sum_{m=1}^{M} (b_{m}(x_j) - \mu(x_j))^2, 
$$
$$
\gamma_1(x_j) = \frac{1}{M} \sum_{m=1}^{M} \bigl(\frac{b_{m}(x_j) - \mu(x_j)}{\sigma}\bigr)^3,\qquad  \gamma_2(x_j) = \frac{1}{M} \sum_{m=1}^{M} \bigl(\frac{b_{m}(x_j) - \mu(x_j)}{\sigma}\bigr)^4 - 3.
$$
Here, $\mu$, $\sigma$, $\gamma_1$, and $\gamma_2$ are sample mean, sample standard deviation, sample skewness, and sample excess kurtosis, respectively. Figures \ref{mean_std} shows the sample first four moments of the parameter $b(x)$ versus $x$, and their histograms.

Obviously, the moments are not constant and vary rapidly in $x$. This shows that the uncertain parameter $b(x)$ cannot be accurately represented by stationary random fields. It also shows that the field is not Gaussian, since for instance the skewness is not zero. 
We therefore need to consider other models than stationary Gaussian random fields to have a better characterization of uncertainty in the parameter $a(x)$. One option is to construct a non-stationary and non-Gaussian random field with variable moments. There are however three problems with this option. One main issue is that due to the imprecise character of uncertainties in the problem, it is not possible to accurately describe the moments by crisp values. Secondly, it is often difficult in practice to work with general non-stationary random fields. For instance, a Karhunen-Loeve approximation of such general fields may require iterative solvers based on density estimation approaches, see e.g. \cite{non_Gaussian1,non_Gaussian2}. Moreover, from a computational point of view, numerical stochastic homogenization techniques, which is necessary to treat the multiscale nature of the parameter, are not applicable in the case of non-stationary fields \cite{Souganidis:99,Bourgeata_Piatnitski:04,Blanc_LeBris:10,Gloria_Otto:11}. 

In Section \ref{sec:modeling}, we propose an alternative option, where instead of considering and dealing with a non-stationary random field with variable crisp moments, a fuzzy-stationary random field is considered, i.e. a random field whose first four moments are fuzzy variables and independent of $x$. Such models are both easy to construct and allow for including imprecise uncertainty in the problem. 


\subsection{Correlation length and correlation function}
\label{sec:corr_length}

\begin{figure}[!h]
\vskip -.3cm
  \begin{center}
       \includegraphics[width=0.65\linewidth]{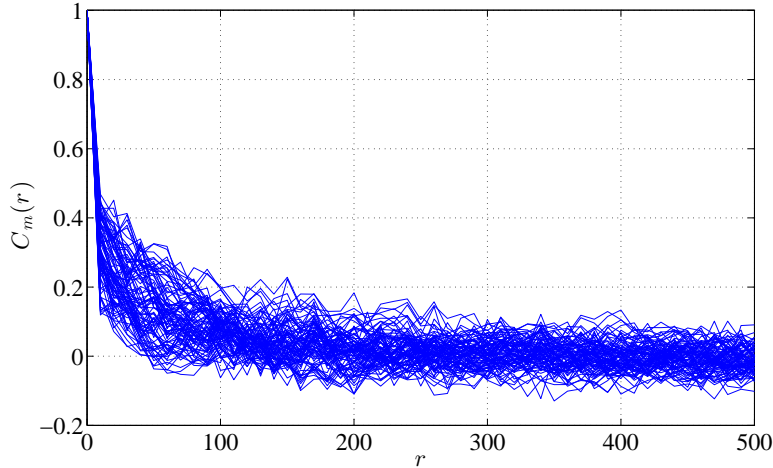}
        \vskip -0.3cm        
\caption{Empirical correlation function $C_m(r_n)$ versus the distance $r_n$ for the parameter $b_m(x)$ for all realizations $m=1,\dotsc, M$.}
\label{CORR1}
  \end{center}
  \vskip -.5cm
\end{figure}

We now study the correlation length $\ell$ of the one-dimensional field $b(x)$. Let $b_{m}(x)$ be the $m$-th realization of the field, corresponding to the $m$-th thin bar, given on a set of $N_x$ discrete points $\{ x_j \}_{j=1}^{N_x}$ along the thin bar. The mean of the field can be estimated by the average of the $m$-th sample over all discrete points
$$
\hat{\mu}_m = \frac{1}{N_x} \sum_{j=1}^{N_x} b_{m}(x_j). 
$$
For an increasing sequence of distances $r_n = n \, h$, with $n=1,2,3,\dotsc$, and $h = 10 \, \mu$m, we find all pairs $(x_i,x_j)$ separated by distance $r_n \pm \frac{h}{2}$: ${\text{find}} \, \, (i,j) \, \, {\text{s.t.}} \, \, | x_i - x_j | \in [ r_n - \frac{h}{2}, r_n + \frac{h}{2} ]$. 
For such pairs, we introduce the normalized empirical correlation function
$$
C_{m}(r_n) = \frac{\sum_{(i,j)} \bigl( b_{m}(x_i) -  \hat{\mu}_m \bigr) \, \bigl( b_{m}(x_j) -  \hat{\mu}_m \bigr)}{\sqrt{ \sum_i \bigl( b_{m}(x_i) -  \hat{\mu}_m \bigr)^2} \, \sqrt{ \sum_j \bigl( b_{m}(x_j) -  \hat{\mu}_m \bigr)^2}}.
$$
We repeat the above procedure for all $M$ realizations. Figure \ref{CORR1} shows the normalized correlation function $C_m(r_n)$ versus $r_n$ for all realizations $m=1, \dotsc, M$. 

We observe that the correlation function decreases from 1 at $r=0$ by $50 \%$ at $r=10$, by $70 \%$ at a distance close to $r = 50$, by $80 \%$ at a distance close to $r = 100$, and continues decreasing by more than $90 \%$ at larger distances. This suggests that the correlation length of the parameter field $b(x)$ is of the order of fiber diameters ($\sim 5 \, h - 20 \, h$). We will use this observation in Section \ref{sec:modeling} to construct the empirical field. In practice, the value of correlation length for the empirical field must be selected so that the best fit to the data is obtained.



\section{Fuzzy set theory}
\label{sec:fuzzy}

Fuzzy sets, introduced by Lotfi A. Zadeh \cite{Zadeh:65}, generalize the notion of classical sets. In a classical set, a single element can either belong or not belong to the set. The notion of membership in this case can be viewed as a characteristic function, taking values either 0 or 1 to characterize the absence or presence of elements in the set. In fuzzy set theory, a grade of membership between 0 and 1 is assigned to each element. The more an element belongs to the fuzzy set, the closer its grade of membership is to 1. Fuzzy sets were initially introduced to accommodate fuzziness contained in human language, judgments, and decisions. Later, it was recognized that the membership degree can represent a degree of possibility \cite{Zadeh:78,DuboisPrade:88}, which allows describing certain forms of uncertainty that is inherently non-statistical.

This section provides a brief overview of the basic mathematical framework of fuzzy set theory. Only the concepts relevant to the focus of this work are mentioned here. For a detailed description, we refer for instance to \cite{Ross:2010,Dubois_Prade:1980,Zimmermann}.

\subsection{One-dimensional fuzzy variables}

Let $Z \subset {\mathbb R}$ be a classical set of elements, called the universe or space, whose generic elements are denoted by $z$. A fuzzy set (or fuzzy variable) is defined by a set $\tilde{z}$ of pairs
$$
\tilde{z} = \{ (z, \mu(z)), \, \, z \in Z \},
$$
where $\mu(z): {\mathbb R} \rightarrow [0,1]$ is a continuous membership function (or generalized characteristic function). We note that the particular case with $\sup_{z \in Z} \mu(z) = 1$ considered here is referred to as a normalized membership function. More general membership functions which are piecewise continuous and non-normalized can be also considered. We further assume that the fuzzy set is convex.

\begin{figure}[!h]
\psfrag{mu}[c][][0.8]{$\mu(z)$}
\psfrag{z}[c][][0.8]{$z$}
  \begin{center}
        \subfigure[convex fuzzy set]{\includegraphics[width=7cm,height=5cm]{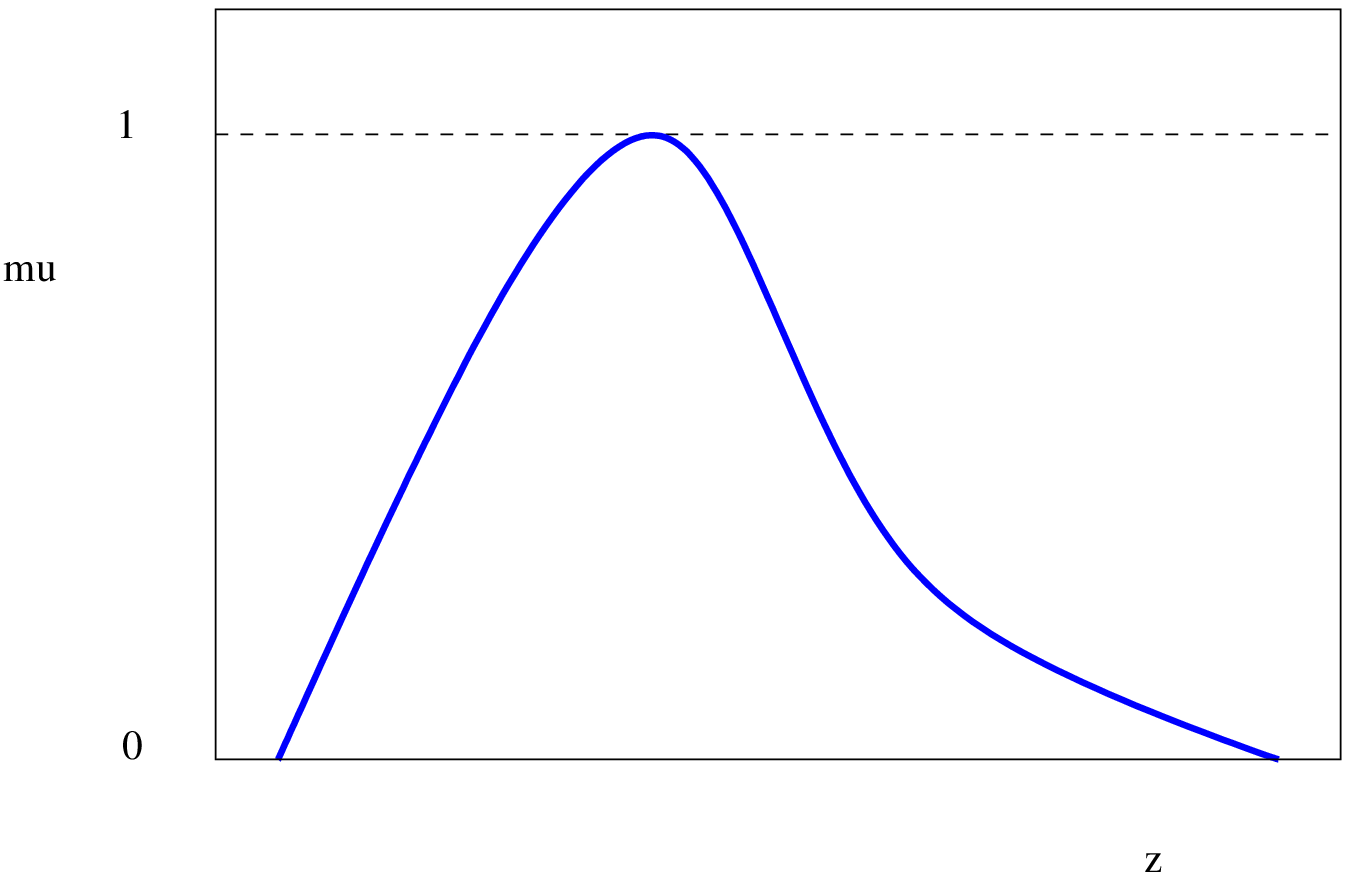}}
        \hskip 1cm
        \subfigure[non-convex fuzzy set]{\includegraphics[width=7cm,height=5cm]{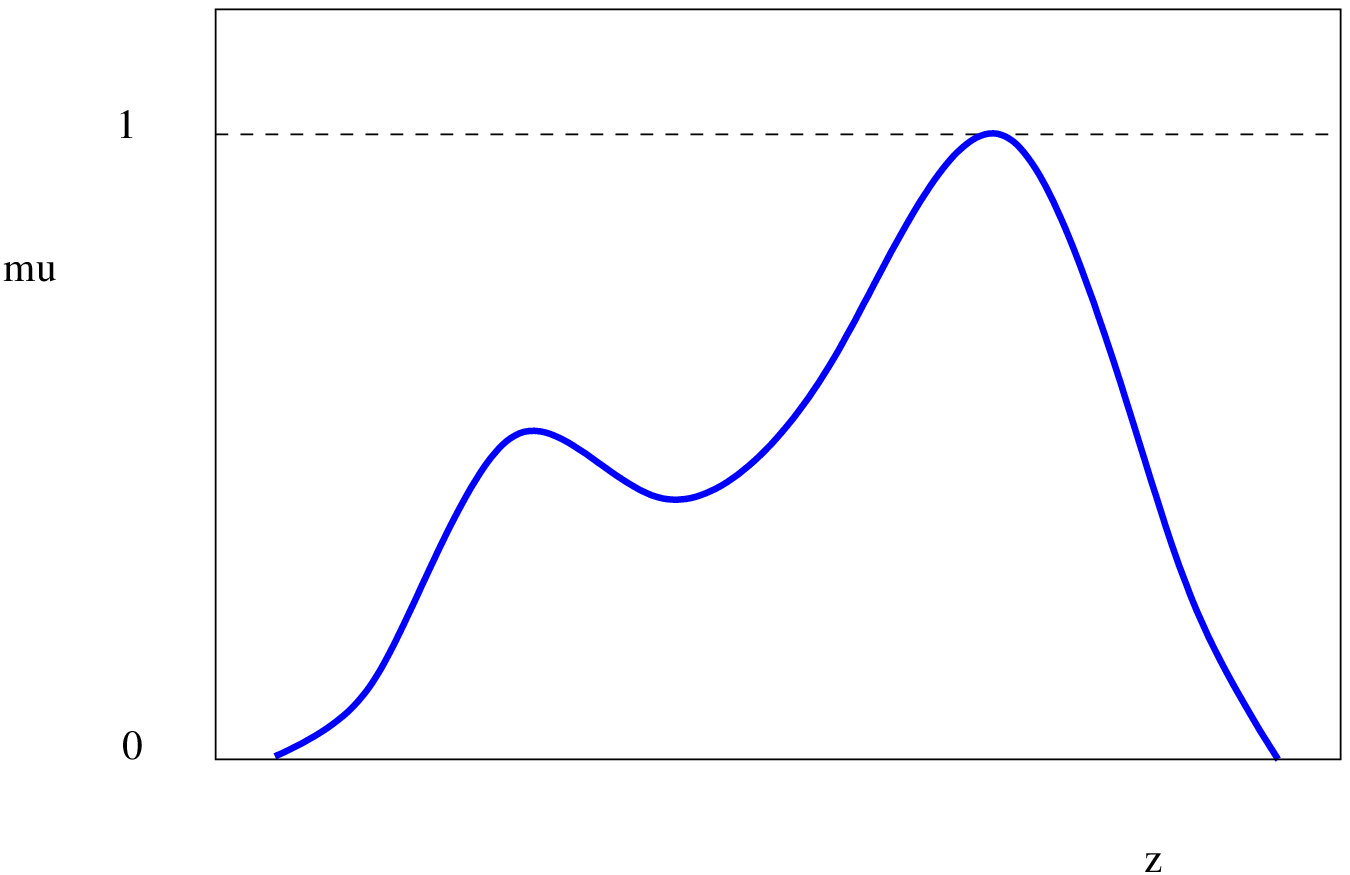}}
\vskip -0.3cm        
\caption{Membership functions of convex (left) and non-convex (right) fuzzy sets.}
\label{convexity}
  \end{center}
  \vskip -.7cm
\end{figure}

\begin{definition}
A fuzzy set $\tilde{z} = \{ (z, \mu(z)), \, \, z \in Z \}$ is convex if 
$$
\forall \, \lambda \in [0,1] \, \, \, \text{and} \, \, \,  \forall \, z_1, z_2 \in Z: \, \, \mu (\lambda \, z_1 + ( 1- \lambda ) \, z_2) \ge \min ( \mu(z_1), \mu(z_2)).
$$
\end{definition}
We note that the convexity of a fuzzy set does not necessarily imply the convexity of the membership function. In fact, according to the definition above, for a fuzzy set to be convex, we need its membership function to be monotonically decreasing on each side of the maximum value. See Figure \ref{convexity}.

\begin{definition}
A fuzzy set $\tilde{z} = \{ (z, \mu(z)), \, \, z \in Z \}$ is positive (resp. negative) if $\mu(z) = 0, \, \forall z <0$ (resp. $\forall z > 0$). This is denoted by $\tilde{z} > 0 $ (resp. $\tilde{z}<0$).
\end{definition}

An important notion in fuzzy set theory is the notion of $\alpha$-cuts, which allows to decompose fuzzy computations into several interval computations. By employing interval arithmetic and computations \cite{Interval:01,Interval:09}, we can introduce mathematical operations for fuzzy variables.

\medskip
\noindent
$\boldsymbol \alpha${\bf -cut representation.} 
Let $\tilde{z} =\{ (z, \mu(z)), \, \, z \in Z \}$ be a convex fuzzy set with a continuous membership function $\mu(z): {\mathbb R} \rightarrow [0,1]$. 
The membership function $\mu(z)$ can be identified with the one-parametric family of sets
$$
{\mathcal S} := \{ S(\alpha) \subset {\mathbb R} \, | \, \alpha \in [0,1] \},
$$
where the $\alpha$-cut or $\alpha$-level set $S(\alpha)$ is defined as 
\begin{equation}\label{alpha-set}
\forall \alpha \in (0,1]: \, \, \,  S(\alpha) =  \{ z \in  Z \, | \, \mu(z) \ge \alpha \}, \quad \text{and} \quad S(0) = {\text{closure}}\{ z \in Z \, | \, \mu(z) >0 \},
\end{equation} 
with the following properties:
\begin{itemize}
\item[1.] $\forall \, \alpha \in [0,1]$, $S(\alpha)$ is a closed set.

\item[2.] $\forall \, \alpha \in [0,1]$, $S(\alpha)$ is a convex (connected) set.

\item[3.] $S(\alpha_1) \subset S(\alpha_2)$ for $\alpha_1 \ge \alpha_2$.
\end{itemize} 

\begin{definition} The support of a fuzzy set $\tilde{z} =\{ (z, \mu(z)), \, \, z \in Z \}$ is the crisp set of all $z \in Z$ such that $\mu(z) > 0$. The closure of the support is the zero-cut $S(0)$.
\end{definition}

\subsection{Multi-dimensional fuzzy variables and their interaction}
\label{sec:multi-dim_fuzzy}

Let $\tilde{\bf z} = \{ ({\bf z}, \mu({\bf z})), \, \, {\bf z} \in {\mathbb R}^n \}$ be a convex $n$-dimensional fuzzy set (or fuzzy vector) with a continuous joint membership function $\mu({\bf z}): {\mathbb R}^n \rightarrow [0,1]$. 
Analogous to the one-dimensional case, the joint membership function $\mu({\bf z})$ can be identified with the one-parametric family of joint $\alpha$-cuts
$$
{\mathcal S} := \{ S(\alpha) \subset {\mathbb R}^n | \alpha \in [0,1] \},
$$
where the joint $\alpha$-cuts are given by
\begin{equation}\label{joint-alpha-set}
\forall \alpha \in (0,1]: \, \, \,  S(\alpha) = \{ {\bf z} \in  {\mathbb R}^n \, | \, \mu({\bf z}) \ge \alpha \}, \quad \text{and} \quad S(0) = {\text{closure}}\{ {\bf z} \in  {\mathbb R}^n \, | \, \mu({\bf z}) >0 \}.
\end{equation}


An important issue that needs to be taken into account when working with two or more fuzzy sets is the interaction (or correlation) between them.
\begin{definition} \label{def_non_interactive} (Non-interactive fuzzy sets) 
Let $\tilde{z}_i = \{ (z_i, \mu_i(z_i)), \, \, z_i \in Z_i \subset {\mathbb R} \}$, with $i = 1, \dotsc, n$, be $n \ge 2$ fuzzy sets. The fuzzy sets are said to be non-interactive if their joint membership function is given by
$$
\mu(z_1, \dotsc, z_n) = \min (\mu_1(z_1), \dotsc, \mu_n(z_n)).
$$
\end{definition}
%

Is it also useful to express non-interactivity of fuzzy sets using the notion of $\alpha$-cuts. The following definition is an important consequence of Definition \ref{def_non_interactive} and is equivalent to it.
\begin{definition} \label{def_non_interactive2}
(Non-interactive fuzzy sets) Let $\tilde{z}_i = \{ (z_i, \mu_i(z_i)), \, \, z_i \in Z_i \subset {\mathbb R} \}$, with $i = 1, \dotsc, n$, be $n \ge 2$ fuzzy sets, and let $S_i(\alpha)$ denote their one-dimensional $\alpha$-cuts, defined by \eqref{alpha-set}. The fuzzy sets are non-interactive if their $n$-dimensional joint $\alpha$-cut corresponding to their joint membership function, defined by \eqref{joint-alpha-set}, is given by the Cartesian product
$$
S(\alpha) = S_1(\alpha) \times \dotsc \times S_n(\alpha), \qquad \forall \alpha \in [0,1],
$$
which is an $n$-dimensional hyperrectangle. 
\end{definition}

Note that in the case when $n=2$, $S(\alpha)$ is a rectangle, and when $n=3$, $S(\alpha)$ is a box. 
%
%
%
%
%
We can further use the notion of $\alpha$-cuts to define interactive fuzzy sets.
\begin{definition} \label{def_interactive}
(Interactive fuzzy sets) Let $\tilde{z}_i = \{ (z_i, \mu_i(z_i)), \, \, z_i \in Z_i \subset {\mathbb R} \}$, with $i = 1, \dotsc, n$, be $n \ge 2$ fuzzy sets, and let $S_i(\alpha)$ denote their one-dimensional $\alpha$-cuts, defined by \eqref{alpha-set}. The fuzzy sets are said to be interactive if their $n$-dimensional joint $\alpha$-cut corresponding to their joint membership function, defined by \eqref{joint-alpha-set}, satisfies
$$
S(\alpha) \subsetneq S_1(\alpha) \times \dotsc \times S_n(\alpha), \qquad \forall \alpha \in [0,1].
$$
Furthermore, if $S(\alpha)$ is the main space diagonal of the hyperrectangle $S_1(\alpha) \times \dotsc \times S_n(\alpha)$, the fuzzy sets are said to be completely interactive. The fuzzy sets are said to be partially interactive if they are neither non-interactive nor completely interactive.
\end{definition}

Note that in the case of complete interaction, when $n=2$, $S(\alpha)$ is the diagonal of a rectangle, and when $n=3$, $S(\alpha)$ is the triagonal of a box. 
We further note that the interaction of fuzzy variables in ${\mathbb R}^n$ is analogous to the correlation of random variables. Moreover, there is not a unique way to represent partial interaction. In the following example, we present a way to express partial interactivity which may be useful for computational purposes. 

\begin{figure}[!h]
  \begin{center}
       \subfigure{\includegraphics[width=5.4cm,height=6cm]{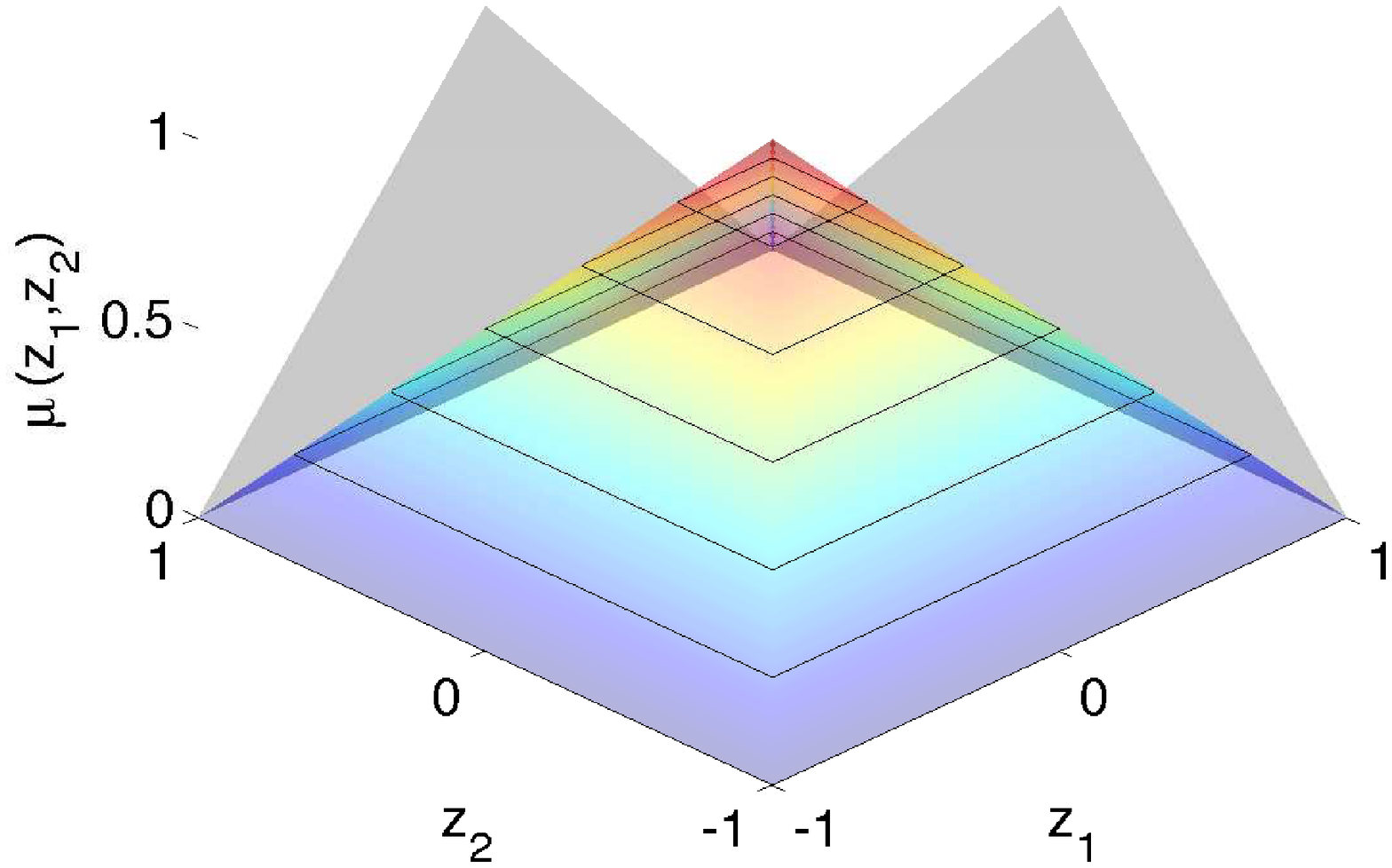}}
        \subfigure{\includegraphics[width=5.4cm,height=6cm]{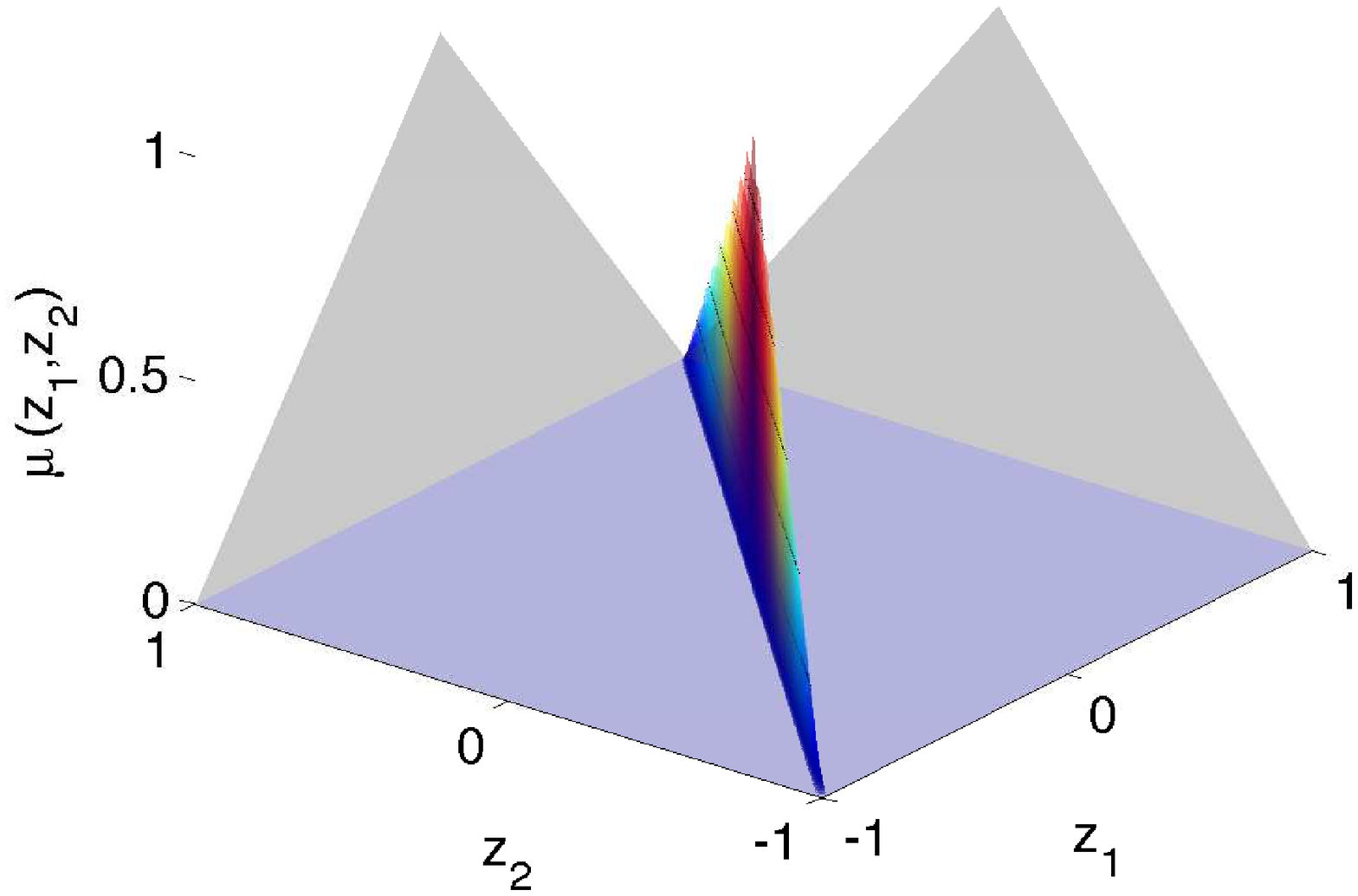}}
        \subfigure{\includegraphics[width=5.4cm,height=6cm]{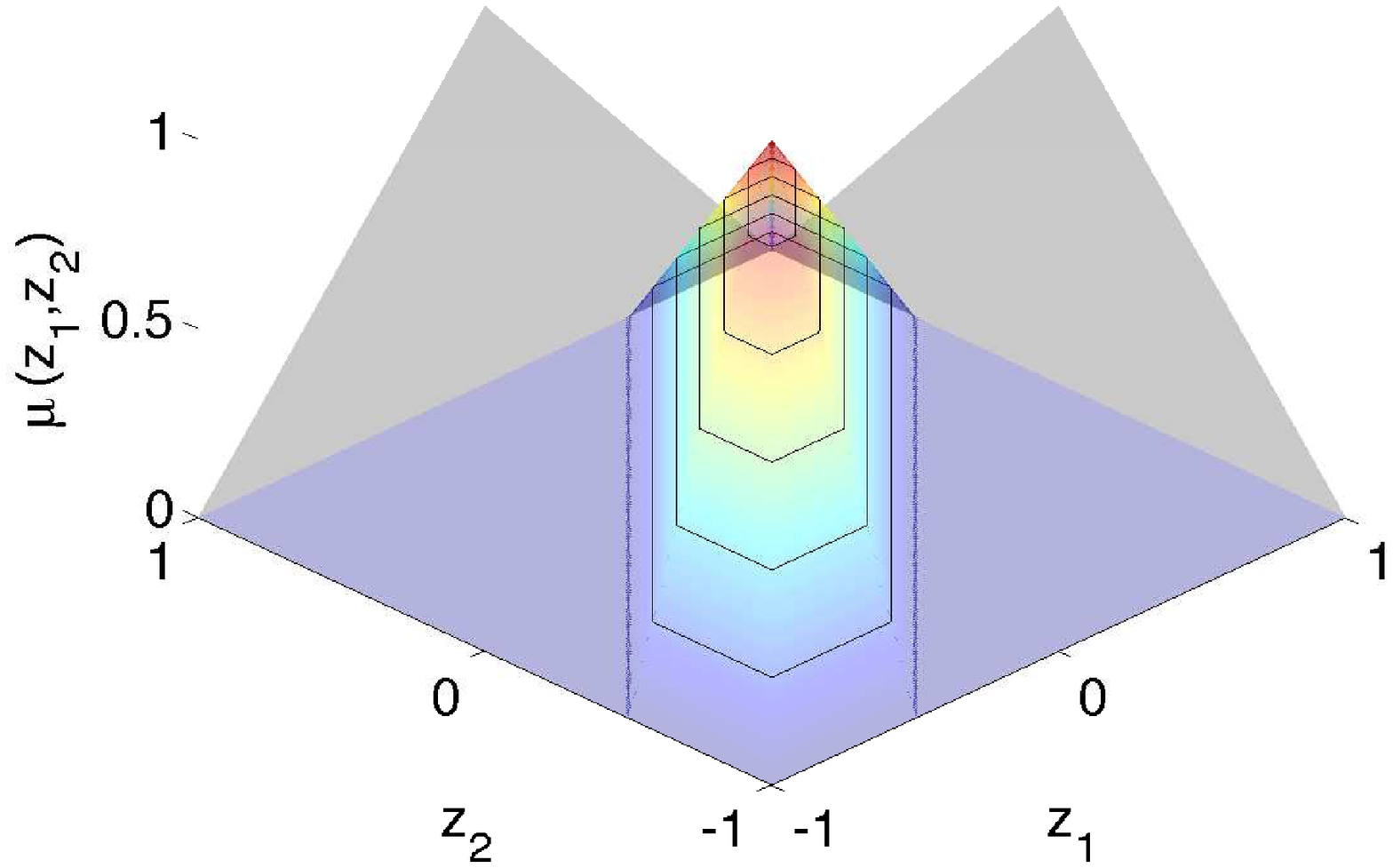}}\\
        \subfigure{\includegraphics[width=5cm,height=3.7cm]{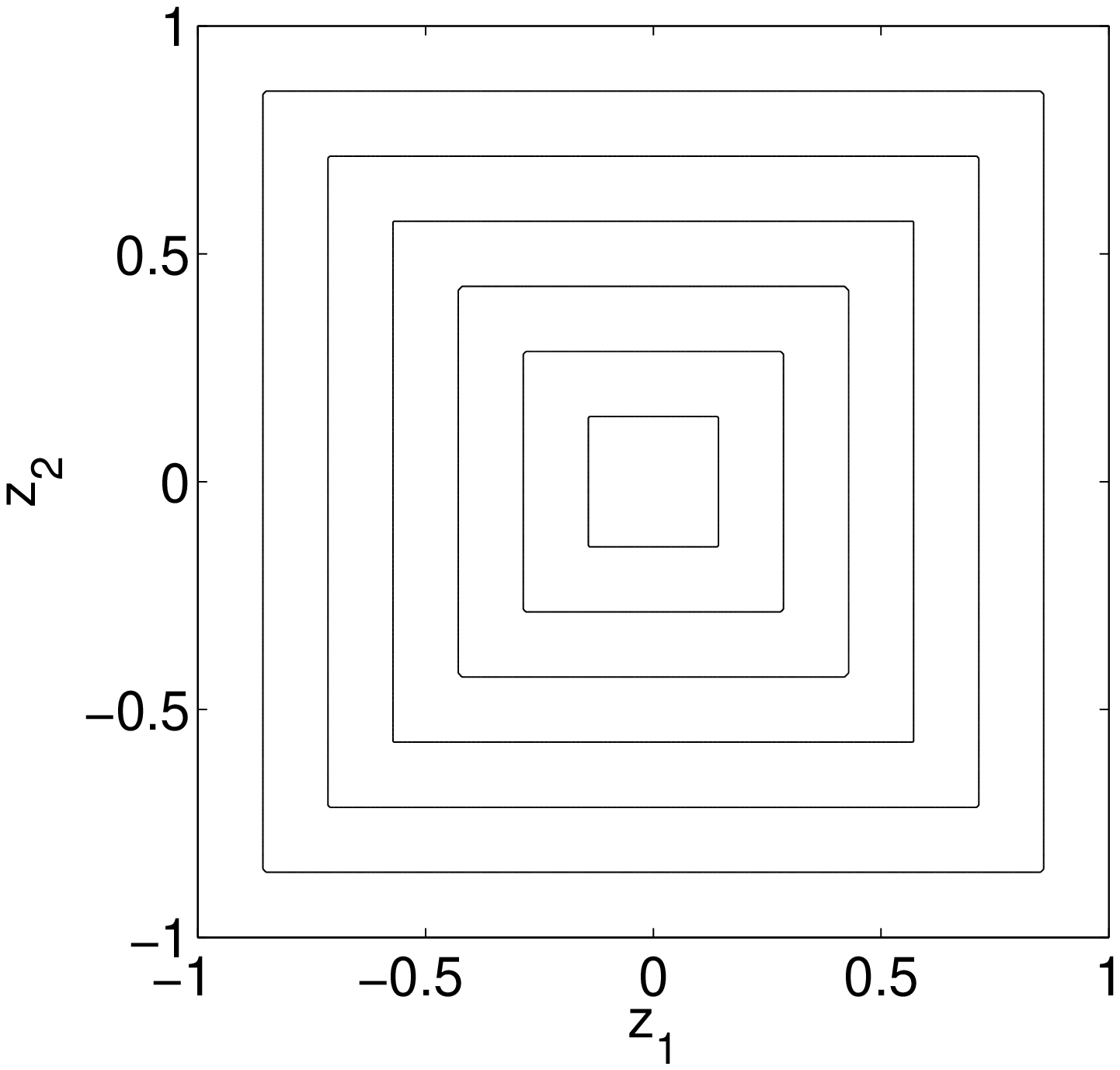}}
        \hskip .5cm
        \subfigure{\includegraphics[width=5cm,height=3.7cm]{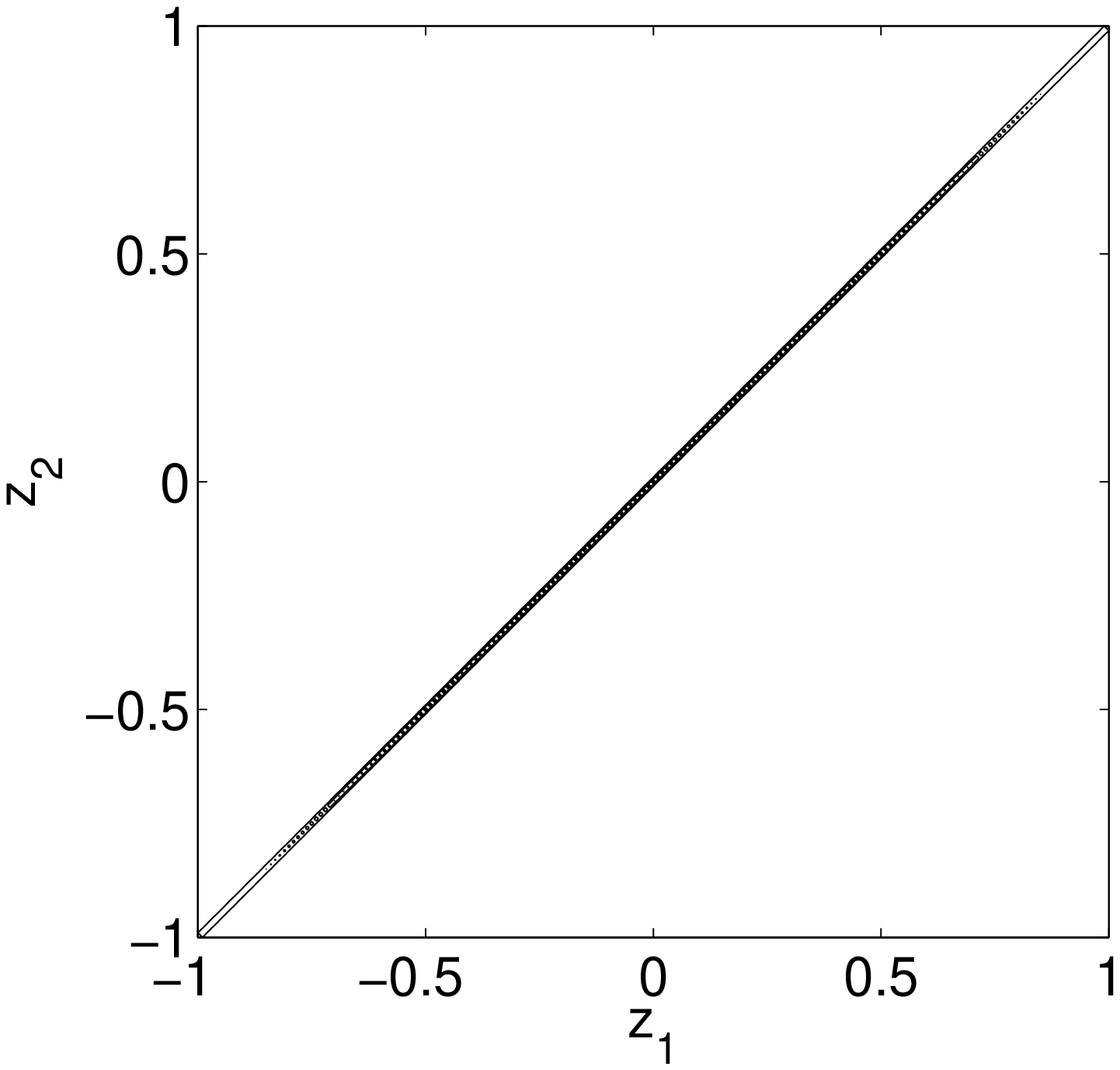}}
        \hskip .5cm
        \subfigure{\includegraphics[width=5cm,height=3.7cm]{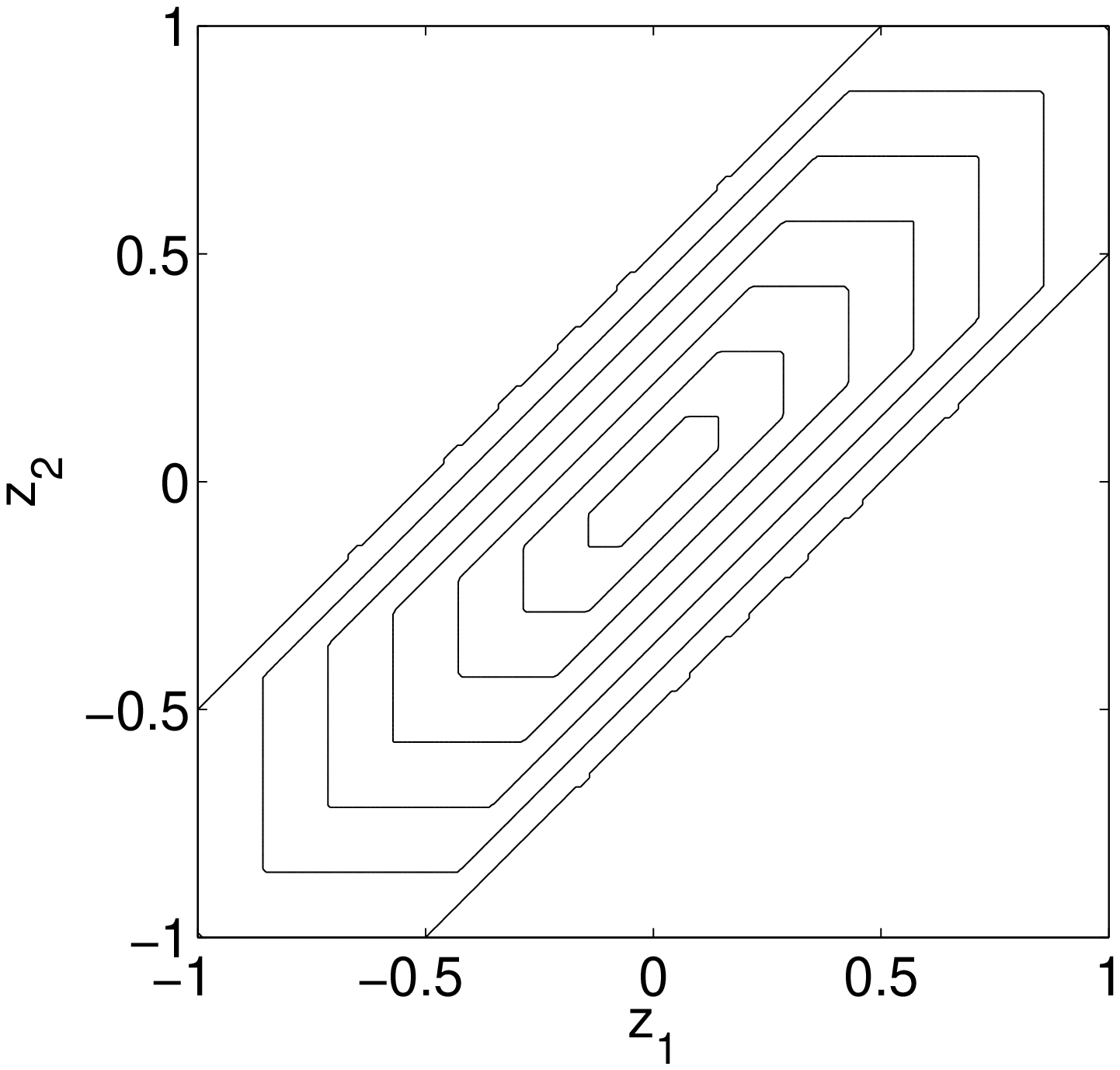}}\\
        \hskip .8cm {\bf Non-interactive} \hskip 1.2cm {\bf Completely interactive} \hskip 1cm {\bf Partially interactive}
\caption{Top figures show the joint membership functions $\mu(z_1,z_2)$, and bottom figures show different $\alpha$-sets $S(\alpha)$, which are the level sets of the membership function, corresponding to different $\alpha$-cuts. The interaction parameter is chosen $\beta=0.5$.}
\label{Interaction_nD}
  \end{center}
  \vskip -.8cm
\end{figure}

\medskip
\noindent
{\bf Example 1.} Let $n=2$ and consider two fuzzy sets $\tilde{z}_i = \{ (z_i, \mu_i(z_i)), \, \, z_i \in Z_i \subset {\mathbb R} \}$, with $i = 1, 2$, with the membership functions
\begin{displaymath}
\mu_1(z_1) =\left\{ \begin{array}{l c}
1- |z_1| & \quad |z_1| \le 1,\\
0 &  \quad |z_1| >1,
\end{array} \right.
\qquad
\mu_2(z_2) =\left\{ \begin{array}{l c}
1- |z_2| & \quad |z_2| \le 1,\\
0 &  \quad |z_2| >1,
\end{array} \right.
\end{displaymath}
We consider three different cases:
\begin{itemize}
\item[{\bf I.}] {\it Non-interactive fuzzy sets}. By Definition \ref{def_non_interactive}, we have
\begin{displaymath}
\mu(z_1,z_2) =\left\{ \begin{array}{l c}
1- \max (|z_1|,|z_2|) & \quad |z_1|, |z_2| \le 1\\
0 &  \quad |z_1|, |z_2| >1
\end{array} \right.
\end{displaymath}
Moreover, we have $S_1(\alpha) = S_2(\alpha)= [\alpha - 1, 1- \alpha]$. Then by Definition \ref{def_non_interactive2}, we have
$$
S(\alpha) = S_1(\alpha) \times S_2(\alpha)= [\alpha - 1, 1- \alpha]^2,
$$
which is a square. See Figure \ref{Interaction_nD} (left column).

\item[{\bf II.}] {\it Completely interactive fuzzy sets}. In the particular case when the membership functions of the two fuzzy sets are the same, as in the present example, we can compute the joint membership function: 
\begin{displaymath}
\mu(z_1,z_2) =\left\{ \begin{array}{l l}
1- |z_1| & \quad z_1 = z_2 \, \, \, \, \, \text{and} \, \, \, \, \, |z_1| \le 1\\
0 &  \quad z_1 \neq z_2 \, \, \, \, \, \text{or} \, \, \, \, \, |z_1|, |z_2| >1
\end{array} \right.
\end{displaymath}
By Definition \ref{def_interactive}, $S(\alpha)$ is the main diagonal of the square $[\alpha - 1, 1- \alpha]^2$, which passes through the origin and is given by equation $z_1 = z_2$. See Figure \ref{Interaction_nD} (middle column).

\item[{\bf III.}] {\it Partially interactive fuzzy sets}. In this case, an $\alpha$-cut needs to be a geometric shape between a square (non-interactive) and its diagonal (completely interactive). For instance, we let $S(\alpha)$ be a hexagon around the diagonal with a parameter $\beta>0$, representing the thickness of the hexagon orthogonal to the diagonal, referred to as interaction parameter. See Figure \ref{Interaction_nD} (right column). The joint membership function in this simple case is given by
\begin{displaymath}
\mu_{\beta}(z_1,z_2) =\left\{ \begin{array}{l c}
\max \Bigl( 0, \min \bigl( 1- \frac{|z_1-z_2|}{\beta}, 1- \max (|z_1|,|z_2|) \bigr) \Bigr)& \quad |z_1|, |z_2| \le 1\\
0 &  \quad |z_1|, |z_2| >1
\end{array} \right.
\end{displaymath}
\end{itemize}
%

\begin{remark}
We note that the hexagon representation proposed above can be generalized to more complicated cases, for example when the two fuzzy sets have different membership functions, or when there are $n \ge 3$ fuzzy sets. One however needs to introduce more interaction parameters $\beta$. For instance when $n=3$ and all three fuzzy sets $\tilde{z}_1,\tilde{z}_2$, and $\tilde{z}_3$ are partially interactive and with the same triangular membership functions, we need to choose three different parameters $\beta_1, \beta_2, \beta_3$. 
If among the three fuzzy sets only two of them are partially interactive, then we need only one interaction parameter.
\end{remark}

\begin{remark}
In the particular case when the membership functions of two completely interactive fuzzy variables are the same, i.e. $\mu_1(z_1) \equiv \mu_2(z_2)$, (like in Example 2, Case II), their joint membership function can be directly computed
\begin{displaymath}
\mu(z_1,z_2) =\left\{ \begin{array}{l l}
\mu(z_1) & \quad z_1 = z_2 \\
0 &  \quad z_1 \neq z_2
\end{array} \right.
\end{displaymath}
In this case, the joint $\alpha$-cut passes through the origin $z_1 = z_2 =0$ and is given by the line equation $z_1 = z_2$. However, if the membership functions are not the same, the direct derivation of the joint membership function is cumbersome. In this case, the joint $\alpha$-cut is still a straight line, but it does not pass through the origin.
\end{remark}

\subsection{Fuzzy functions and fuzzy operators}
\label{sec:fuzzy_vec}

A fuzzy function is a generalization of the concept of a classical function. A classical function is a mapping from its domain of definition into its range. Two important generalizations include: 1) a crisp map with fuzzy arguments; and 2) a fuzzy map with crisp arguments. Both cases generate an output fuzzy set. In the present work, we only consider the first case and refer to a crisp map with fuzzy arguments as a {\it fuzzy function}.
%
%
The definition of fuzzy functions is based on a fundamental axiom in fuzzy set theory, known as {\it generalized extension principle}. It provides a general method for extending non-fuzzy mathematical concepts in order to deal with fuzzy quantities. 
\begin{definition} \label{generalized_extension} (Generalized extension principle \cite{Fuller:04,Fuller:14}) Let $\tilde{\bf z} = \{ ({\bf z}, \mu({\bf z})), \, \, {\bf z} \in Z \subset {\mathbb R}^n \}$ be a convex $n$-dimensional fuzzy set with a continuous joint membership function $\mu({\bf z}): Z \rightarrow [0,1]$. Let $g: Z \rightarrow W$ be a continuous mapping such that ${\mathsf w} = g({\bf z})$. The membership function of the output fuzzy set ${\mathsf w}$ in space $W$ is given by 
\begin{displaymath}
\mu_{{\mathsf w}}({\mathsf w}) =\left\{ \begin{array}{l l}
\sup_{{\bf z}= g^{-1}({\mathsf w})}  \, \mu({\bf z}) & \qquad g^{-1}({\mathsf w}) \neq \emptyset \\
0 &  \qquad g^{-1}({\mathsf w}) = \emptyset
\end{array} \right.
\end{displaymath}
where $g^{-1}({\mathsf w})$ is the inverse image of ${\mathsf w}$, and $\emptyset$ is the empty set.
\end{definition}

We note that the generalized extension principle uses the concept of joint membership function and hence is valid for both non-interactive and interactive fuzzy variables. In the case of non-interactive fuzzy variables, the generalized extension principle turns into Zadeh's sup-min extension principle \cite{Zadeh:75}.

The direct application of the extension principle to perform operations on fuzzy functions can be quite complicated and numerically cumbersome, particularly in the case when $n \ge 2$, due to the interaction between fuzzy variables. The computations can be simplified using the $\alpha$-cut representation of joint membership functions, thanks to the following important result.
\begin{theorem} \label{thm1} ($\alpha$-cut representation of the extension principle \cite{Fuller:04,Fuller:14}) 
Let $\tilde{\bf z} = \{ ({\bf z}, \mu({\bf z})), \, \, {\bf z} \in Z \subset {\mathbb R}^n \}$ be a convex $n$-dimensional fuzzy set with a continuous joint membership function $\mu({\bf z}): Z \rightarrow [0,1]$ and its corresponding joint $\alpha$-cut $S(\alpha)$. Let $g: Z \rightarrow W$ be a continuous mapping such that ${\mathsf w} = g({\bf z})$. Then the $\alpha$-cut $S_{{\mathsf w}}(\alpha)$ corresponding to the membership function of the output fuzzy variable is given by 
$$
S_{{\mathsf w}}(\alpha) = g (S(\alpha)), \qquad \forall \alpha \in [0,1].
$$
\end{theorem}
Motivated by Theorem \ref{thm1}, we present a practical approach for computing fuzzy functions and fuzzy operations based on the $\alpha$-cut representation of the joint membership function and using two main tools: 1) the $\alpha$-cut-based interaction strategy discussed in Section \ref{sec:multi-dim_fuzzy}; and 2) the worst case scenario. We take two major steps to find the output joint $\alpha$-cut:
\begin{itemize}

\item[1.] {\it Interaction}: find the joint $\alpha$-cut $S(\alpha) \subseteq S_1(\alpha) \times \dotsc S_n(\alpha)$ corresponding to the joint membership function of the $n$ input fuzzy variables based on their interaction, where $S_i(\alpha)$, with $i=1, \dotsc, n$, is the one-dimensional $\alpha$-cut of the input fuzzy set $\tilde{z}_i$.

\item[2.] {\it Worst case scenario}: find the $\alpha$-cut $S_{{\mathsf w}}(\alpha)$ corresponding to the membership function $\mu_{{\mathsf w}}({\mathsf w})$ of the output fuzzy variable based on the worst case scenario of all possible outcomes ${\mathsf w}^{(j)} = g({\bf z}^{(j)})$ at all points ${\bf z}^{(j)} \in S(\alpha)$ in the joint $\alpha$-cut of input fuzzy sets.
\end{itemize}
Using the proposed approach above, we can compute fuzzy functions and fuzzy operations. In particular, we define the four basic arithmetic operations on fuzzy functions.

\medskip
\noindent
{\bf Fuzzy arithmetic operations.} 
Consider $n$ fuzzy sets $\tilde{z}_i = \{ (z_i, \mu_i(z_i)), \, \, z_i \in Z_i \subset {\mathbb R} \}$, with $i = 1, \dotsc, n$, and assume that their joint $\alpha$-cut $S(\alpha)$ is computed based on their interaction. Let $g_1$ and $g_2$ be two continuous fuzzy functions with the fuzzy input vector ${\bf z} = (z_1, \dotsc, z_n)$. We want to find the $\alpha$-cuts of the following fuzzy quantities:
$$
{\mathsf w}_1 = g_1({\bf z}) + g_2({\bf z}), \quad {\mathsf w}_2 = g_1({\bf z}) - g_2({\bf z}), \quad {\mathsf w}_3 = g_1({\bf z}) \, g_2({\bf z}), \quad {\mathsf w}_4 = g_1({\bf z}) / g_2({\bf z}).
$$
We follow the proposed approach and first consider all points ${\bf z}_{\alpha}^{(j)} = (z_{1 \alpha}^{(j)}, \dotsc, z_{n \alpha}^{(j)}) \in S(\alpha)$ in the joint $\alpha$-cut of the input fuzzy sets, which is found based on the given interaction of fuzzy sets. Next we compute
$$
\ubar{g_1} :=  \min_{{\bf z}_{\alpha}^{(j)} \in S(\alpha)} g_1({\bf z}_{\alpha}^{(j)}), \, \, \, \, \,
\bar{g_1} :=  \max_{{\bf z}_{\alpha}^{(j)} \in S(\alpha)} g_1({\bf z}_{\alpha}^{(j)}), \, \, \, \, \,
\ubar{g_2} :=\min_{{\bf z}_{\alpha}^{(j)} \in S(\alpha)} g_2({\bf z}_{\alpha}^{(j)}), \, \, \, \, \,
\bar{g_2} :=  \max_{{\bf z}_{\alpha}^{(j)} \in S(\alpha)} g_2({\bf z}_{\alpha}^{(j)}).
$$
Finally, based on the worst case scenario, the output $\alpha$-cuts are given as follows:
%
\begin{align*}
&S_{{\mathsf w}_1}(\alpha) = \bigl[  \ubar{g_1} + \ubar{g_2}, \,  \bar{g_1} + \bar{g_2}    \bigr],\\
&S_{{\mathsf w}_2}(\alpha) = \bigl[  \ubar{g_1} - \bar{g_2}, \,  \bar{g_1} - \ubar{g_2}    \bigr],\\
&S_{{\mathsf w}_3}(\alpha) = \bigl[  \min \bigl( \ubar{g_1} \, \ubar{g_2}, \ubar{g_1} \, \bar{g_2},\bar{g_1} \, \ubar{g_2}, \bar{g_1} \, \bar{g_2}     \bigr), \,  \max \bigl( \ubar{g_1} \, \ubar{g_2}, \ubar{g_1} \, \bar{g_2},\bar{g_1} \, \ubar{g_2}, \bar{g_1} \, \bar{g_2}     \bigr)   \bigr],\\
&S_{{\mathsf w}_4}(\alpha) = \bigl[  \min \bigl( \ubar{g_1} / \ubar{g_2}, \ubar{g_1} / \bar{g_2},\bar{g_1} / \ubar{g_2}, \bar{g_1} / \bar{g_2}     \bigr), \,  \max \bigl( \ubar{g_1} / \ubar{g_2}, \ubar{g_1} / \bar{g_2},\bar{g_1} / \ubar{g_2}, \bar{g_1} / \bar{g_2}     \bigr)   \bigr].
\end{align*}
We note that the devision operation requires that $0 \notin S_{g_2}(\alpha)$ for all $\alpha \in [0,1]$. An interesting consequence of the arithmetic operators defined above is that if $g_1 \equiv g_2$, then ${\mathsf w}_2 = g_1 - g_1 \neq 0$, ${\mathsf w}_3 = g_1 \, g_1$ is not necessarily positive, and ${\mathsf w}_4 = g_1 / g_1 \neq 1$. See also Cahpter 7 of \cite{Klir:06}.


\medskip
\noindent
{\bf Fuzzy integration.}
The proposed approach can easily be extended to performing other operations on fuzzy functions of $n \ge1$ fuzzy sets. For instance let $\tilde{\bf z} = (\tilde{z}_1,  \dotsc,\tilde{z}_n)$ be a fuzzy vector containing $n$ fuzzy sets with a known joint $\alpha$-cut. Suppose that we want to compute the integral of a fuzzy-valued function $g$ over a crisp interval $[0,1]$:
$$
I = \int_0^1 g(x,{\bf z}) \, dx \approx \Delta x \, \sum_{j=1}^m g(x_j,{\bf z}), \qquad \Delta x = \frac1{m-1}.
$$
We first use the trapezoidal rule (or other quadrature rules) and approximate the crisp integral by a sum of $m$ fuzzy functions $g(x_1, {\bf z}), \dotsc, g(x_m, {\bf z})$. Next, by the definition of addition of fuzzy functions, presented above, we can compute different $\alpha$-cuts of the fuzzy set $I$.


\section{A fuzzy-stochastic model}
\label{sec:modeling}


In this section, we construct and validate a new model---as an alternative to current models based on precise probability---for characterizing uncertainty in the input parameter $a(x)$ by combining stochastic fields and fuzzy variables. Motivated by the solution formula \eqref{exact_sol}, we consider and study the reciprocal of the modulus of elasticity, i.e. $b(x)= a^{-1}(x)$. Our goal is to obtain an empirical model for the field $b(x)$ based on the real data collected in Section \ref{sec:data_collection}. We will motivate and explain the construction of the fuzzy-stochastic model, consisting of a non-Gaussian random field whose first four statistical moments, i.e. mean, variance, skewness, and excess kurtosis are fuzzy variables. 
%
%
%

Throughout this section, we consider a local small domain of length $L=10^4 [\mu \text{m}]$, and construct the model for this domain. In a similar way, models for domains of any size can be constructed. This particular choice will be used in Section \ref{sec:global-local} in the global-local approach.

\subsection{Fuzzy statistical moments}
\label{sec:fuzzy_moments}

In Section \ref{sec:moments}, we studied the first four statistical moments of $b(x)$ on a small domain of length $L=10^4 [\mu \text{m}]$. The data used in the statistical analysis were generated from the original set of data (50 samples of length $1700 [\mu \text{m}]$) by bootstrapping. We generated $M=100$ bootstrap samples of length $L=10^4 [\mu \text{m}]$, each at $N_x=10001$ discrete points $\{ x_j \}_{j=1}^{N_x}$ on the interval [0,L]. We then used the bootstrap samples and obtained the histogram of the four sample moments; see Figure \ref{mean_std}. We now want to model the four moments $\mu(x)$, $\sigma(x)$, $\gamma_1(x)$, and $\gamma_2(x)$, which vary in $x$, by four fuzzy variables $\tilde{z}_1,\tilde{z}_2,\tilde{z}_3,\tilde{z}_4$, respectively. This procedure is called the {\it fuzzification} of moments. 
We will fuzzify the moments using their histograms and construct four membership functions, denoted by $\mu_1(z_1), \dotsc, \mu_4(z_4)$.

\begin{figure}[!h]
  \begin{center}
        \subfigure{\includegraphics[width=8cm,height=3.8cm]{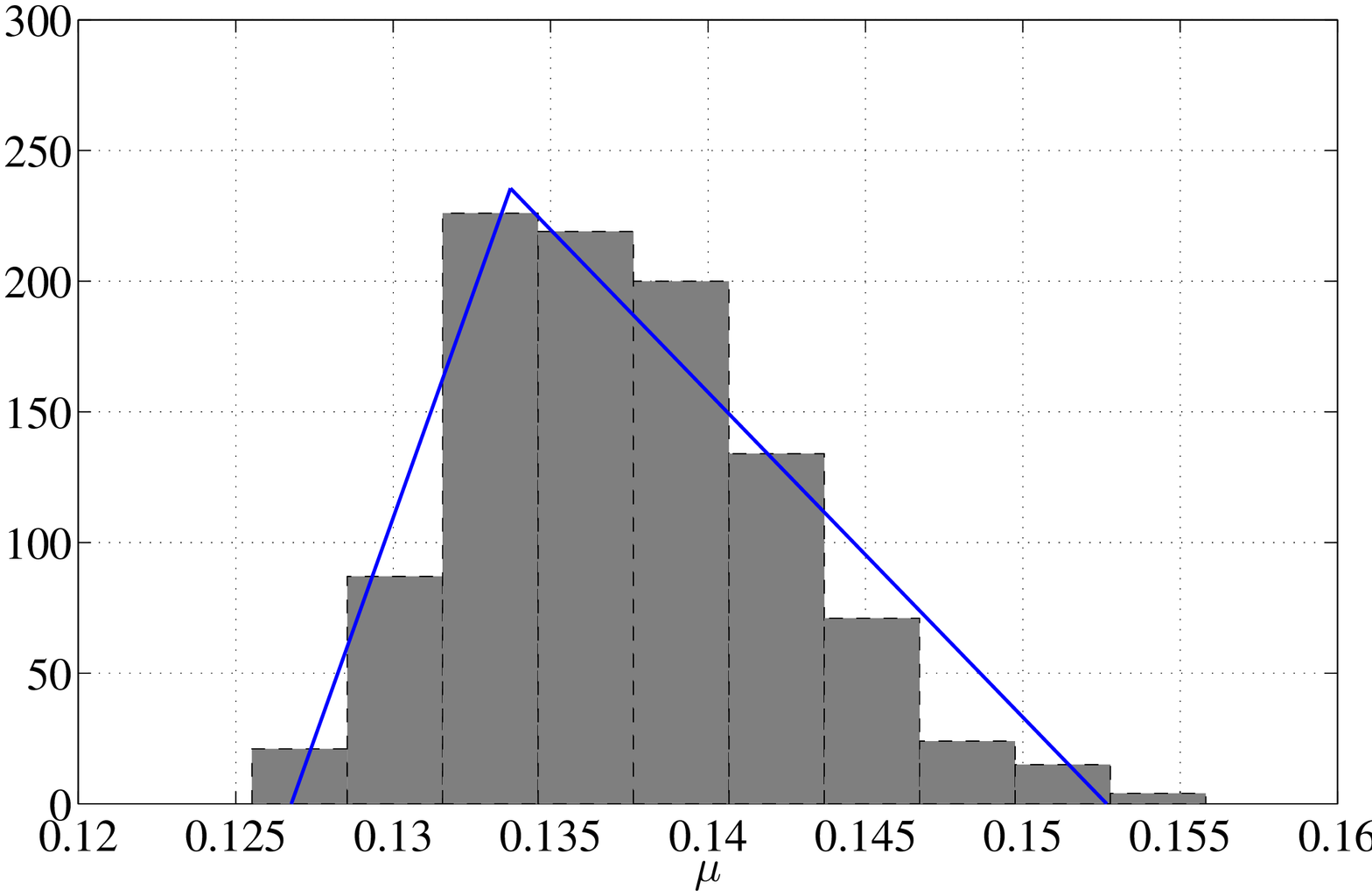}}
        \subfigure{\includegraphics[width=8cm,height=3.8cm]{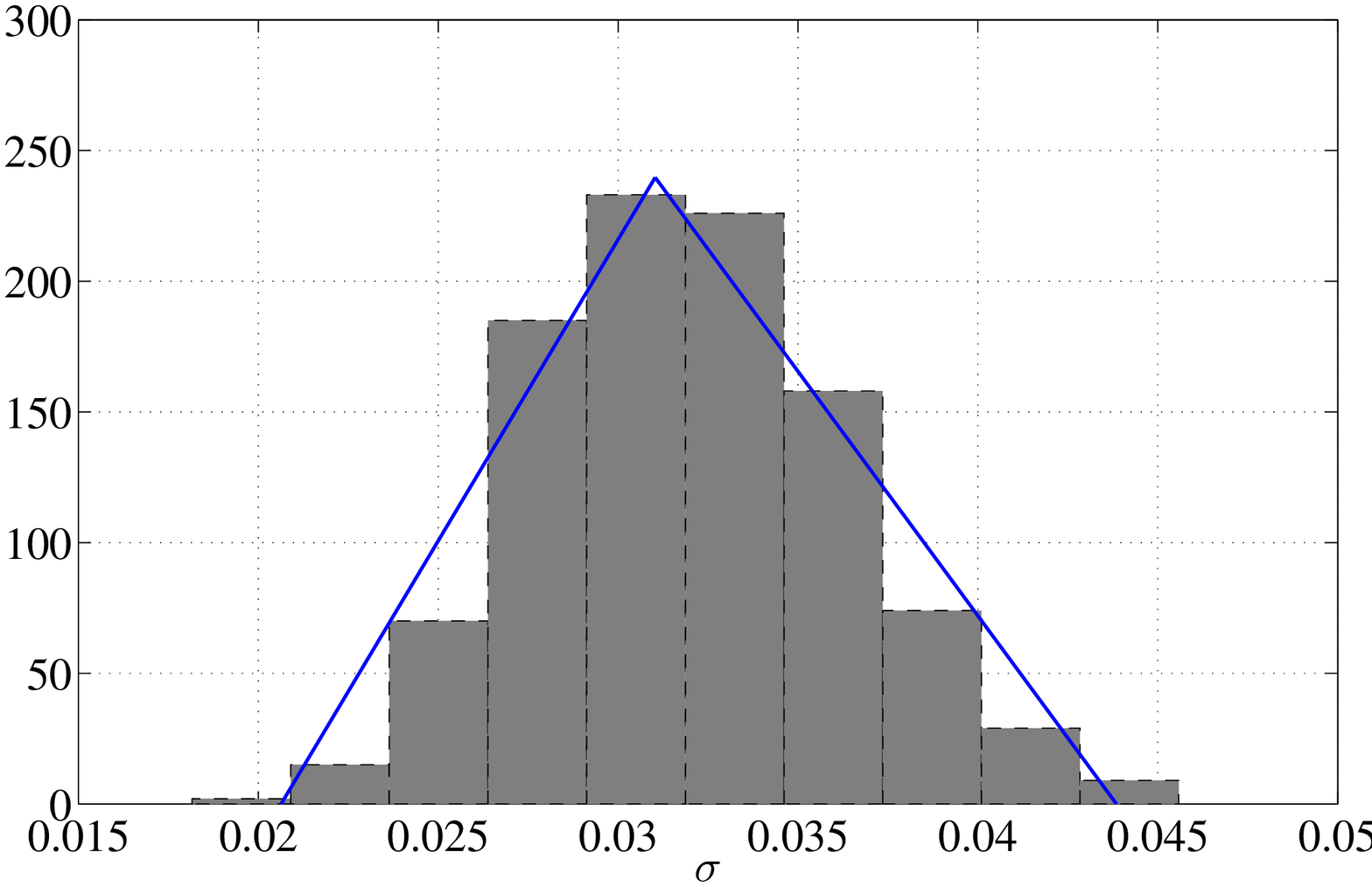}}\\
        \subfigure{\includegraphics[width=8cm,height=3.8cm]{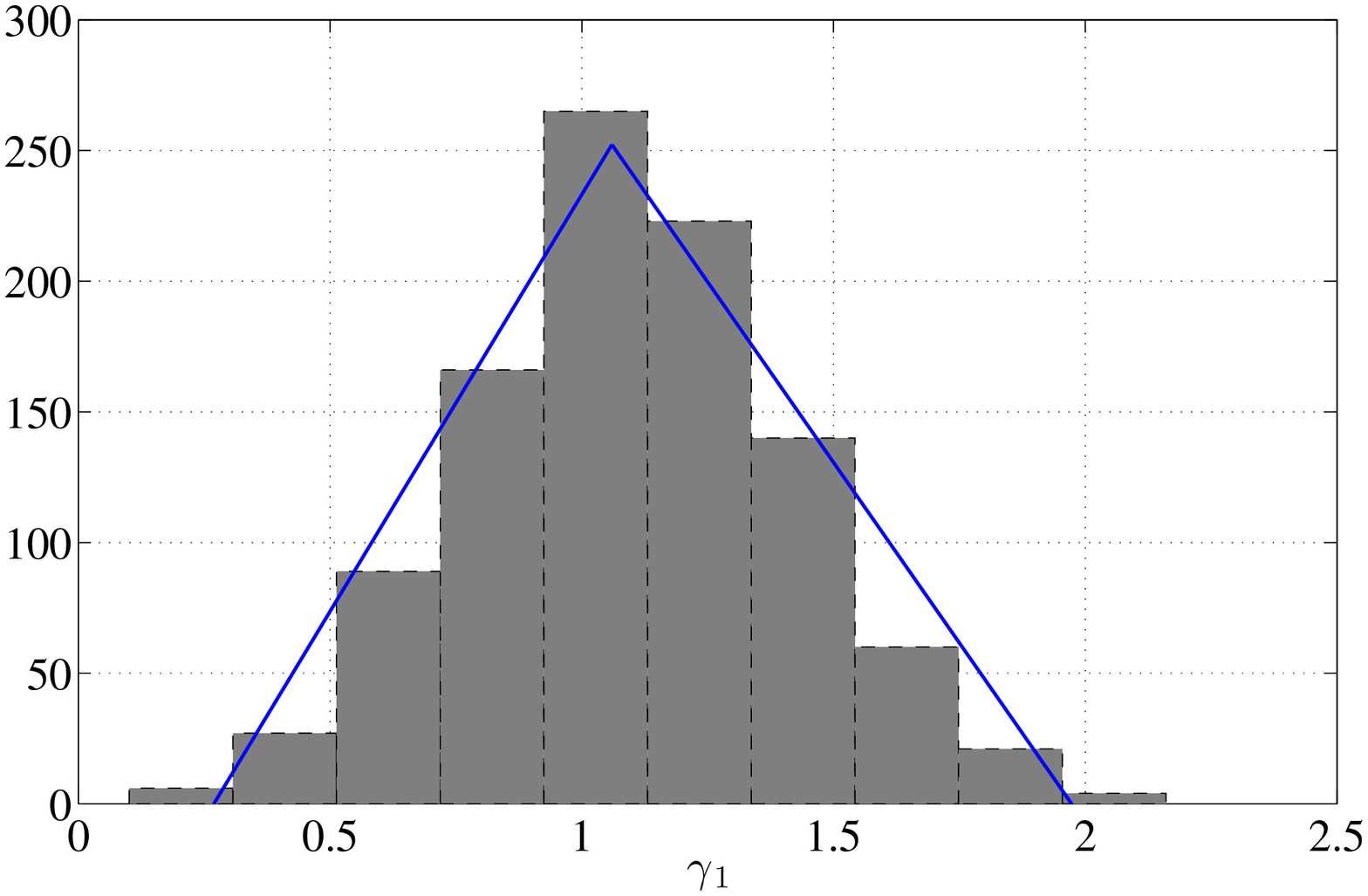}}
        \subfigure{\includegraphics[width=8cm,height=3.8cm]{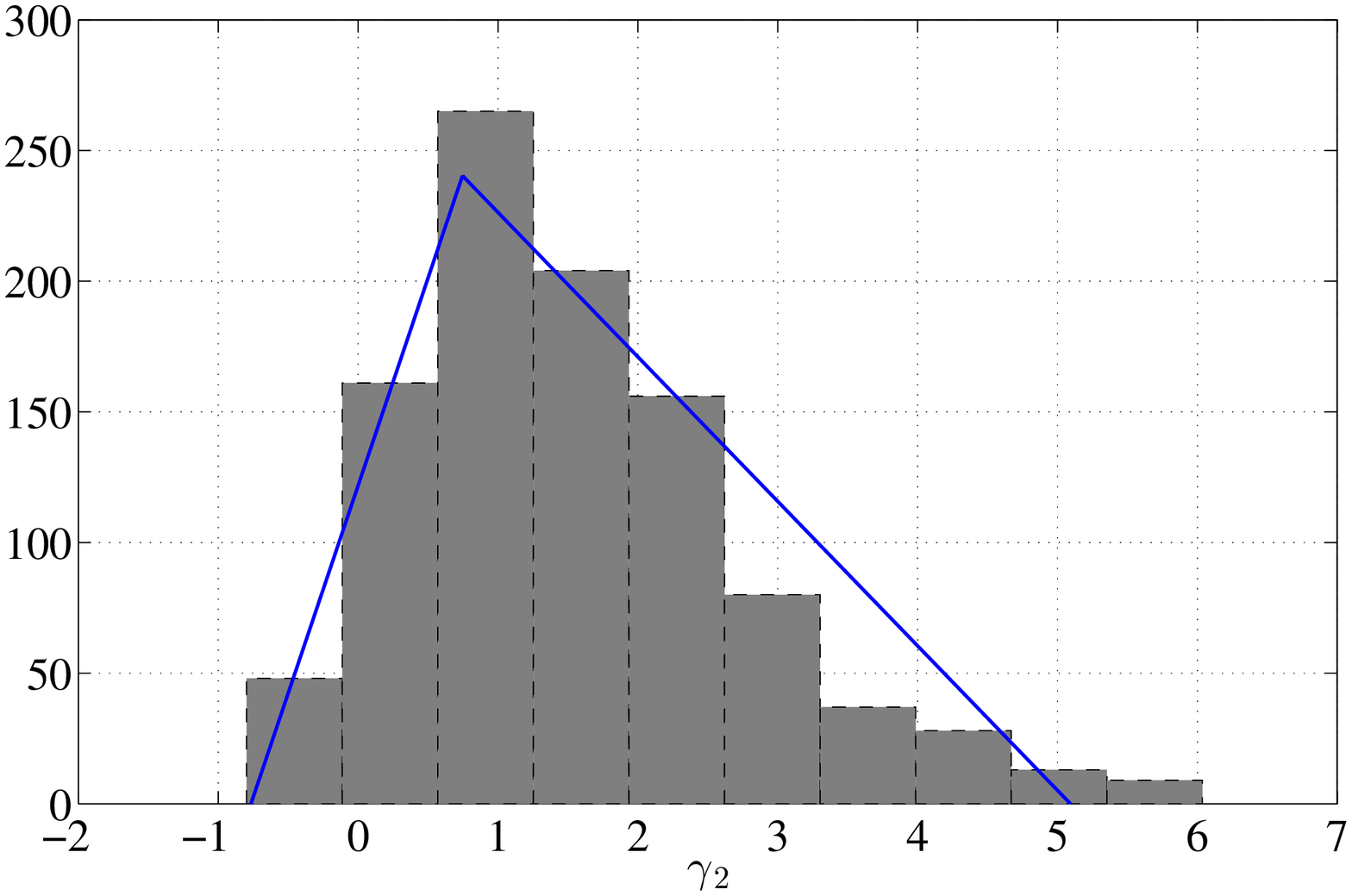}}
\vskip -0.3cm        
\caption{Histograms (gray bars) and non-normalized membership functions (blue solid lines) of sample mean (top left), sample standard deviation (top right), sample skewness (bottom left), and sample excess kurtosis (bottom right) of the field $b(x) = a^{-1}(x)$.} 
\label{hist_membership}
  \end{center}
  \vskip -.3cm
\end{figure}

We first build up the histograms of the moments by dividing the range of the four sample moments (obtained by bootstrap samples) into 10 subintervals and determining the number of sample elements belonging to each subinterval; see Figures \ref{mean_std} and \ref{hist_membership}. As an initial draft for membership functions, we follow  \cite{Moller_Beer:2004} and consider triangular membership functions consisting of two linear branches: a left and a right branch. The two linear branches are determined by the method of least squares as follows. The subinterval containing the largest number of sample elements is selected and called ``mid-subinterval''. The mid-subinterval and all subintervals to its left are used for computing the left branch of the membership function, while the mid-subinterval and all subintervals to its right are used for computing the right branch. We collect two sets of left and right data points, whose $x$- and $y$-coordinates are the mid points of selected subintervals and the corresponding number of sample elements, respectively. We then fit two lines to these two sets of data points by the least squares approach. The intersection of these two lines represent the mean value of the fuzzy variable. Note that the intersecting point may not lie in the mid-interval. The two neighboring zeros of these linear functions mark the interval bounds of the support. See Figure \ref{hist_membership}. It is very important to note that the above procedure is to generate an initial draft for membership functions. We may need to conduct a subsequent modification by imposing additional constraints and corrections. For instance, in Figure \ref{hist_membership} (top left), we may need to change the initial draft of the membership function for the fuzzy set $\tilde{z}_1$ and consider a trapezoidal membership function instead of a triangular function to account for the flat shape of histogram for $\mu \in [0.132, 0.141]$. 
%
%
%
Finally, the (possibly modified) membership functions are normalized so that the function value at the mean value point is one. 



Furthermore, since the four statistical moments are obtained from the same set of data and are directly related to each other, i.e. higher moments are obtained from lower moments, it is reasonable to assume that the four fuzzy variables $\tilde{z}_1,\tilde{z}_2,\tilde{z}_3,\tilde{z}_4$ are completely interactive.

\subsection{A fuzzy-stationary stochastic field}
\label{sec:local-model}
We are now ready to construct an empirical fuzzy-stochastic model for the parameter $b(x)$. 
We denote the empirical fuzzy-stochastic field by $\tilde{b}(x,{\bf y},{\bf z})$, where 
${\bf y} = (y_1, \dotsc,y_N)$ is a random vector containing $N$ random variables, and 
${\bf z} = (z_1,z_2,z_3,z_4)$ is the vector of four fuzzy variables, i.e. the four fuzzy moments obtained in Section \ref{sec:fuzzy_moments}. 

We first note that by the construction, the fuzzy mean of the field is independent of $x$, i.e. ${\mathbb E}[\tilde{b}(x,.,{\bf z})] = z_1$. The other three higher moments are also fuzzy and independent of $x$. We may therefore call the field $\tilde{b}$ a {\it fuzzy-stationary random field}. Moreover, we notice that the field $\tilde{b}$ is not Gaussian, since for instance its skewness is not zero. Motivated by these observations and following the construction of stochastic translation processes \cite{non_Gaussian_stationary}, we propose a fuzzy-stationary non-Gaussian translation random field
\begin{equation} \label{emp_field}
\tilde{b}(x,{\bf y},{\bf z}) = \Psi^{-1}({\bf z}) \circ \Phi(G(x,{\bf y})),
\end{equation}
where $G(x,{\bf y})$ is a stationary Gaussian random field with mean zero and variance one, $\Phi$ is the standard Gaussian cummulative distribution function (CDF), and $\Psi$ is a beta CDF with four parameters, chosen from the four fuzzy moments. 
For the stationary Gaussian field $G(x,{\bf y})$, we consider a squared exponential covariance function: 
\begin{equation}\label{cov_function}
C_G(x_1,x_2) = \exp{\bigl(\frac{-|x_1 - x_2|^2}{2 \, \ell^2}\bigr)}.
\end{equation}
Motivated by the results of Section \ref{sec:corr_length}, we let the filed's spatial correlation length be $\ell = 100$; see Figure \ref{CORR1}. Since the covariance function \eqref{cov_function} is deterministic, we employ the truncated Karhunen-Lo{\'e}ve expansion to represent the stationary Gaussian field:
\begin{equation}\label{KL_N}
G(x,{\bf y}) \approx \sum_{n=1}^{N} \sqrt{\lambda_n} \, \phi_n(x) \, y_n,
\end{equation}
where $\{(\lambda_n, \phi_n(x) )\}_{n=1}^N$ are the eigenpairs of the deterministic covariance function \eqref{cov_function}, and ${\bf y} = (y_1, \dotsc, y_N)  \in {\mathbb R}^N$ is a vector of $N$ uncorrelated normal random variables with zero mean and unit variance, i.e. $y_n \sim {\mathcal N}(0,1)$. The number of terms $N$ is chosen such that a high percentage of the standard deviation is preserved. 
%
%
The recipe for constructing the fuzzy-stochastic field is outlined in Algorithm \ref{ALG_field}.
\begin{algorithm}
\caption{Construction of fuzzy-stochastic field based on the real data} 
\label{ALG_field}
\begin{algorithmic} 
\medskip
\STATE {\bf 0.} Given a length $L$, obtain $M$ bootstrap samples from real data, as in Section \ref{sec:data_collection}.

\medskip
\STATE {\bf 1.} Model first four moments by four fuzzy variables ${\bf z}=(z_1,z_2,z_3,z_4)$, as in Setion \ref{sec:fuzzy_moments}. 

\medskip
\STATE {\bf 2.} Compute a fuzzy beta distribution $\Psi({\bf z})$ with four parameters given by four moments.

\medskip
\STATE {\bf 3.} Generate a stationary Gaussian random field with covariance function \eqref{cov_function} and correlation length $\ell=100$ by the truncated KL expansion \eqref{KL_N}.

\medskip
\STATE {\bf 4.} Obtain the empirical fuzzy-stochastic field by \eqref{emp_field}.
\end{algorithmic}
\end{algorithm}





\noindent
{\bf Model justification.} We now perform a simple justification of the proposed fuzzy-stochastic model by showing that the CDFs of the true field $b(x)$ and the empirical fuzzy-stochastic field $\tilde{b}(x,{\bf y},{\bf z})$ are in good agreement. 
%
%

\begin{figure}[!h]
  \begin{center}
  \includegraphics[width=0.65\textwidth]{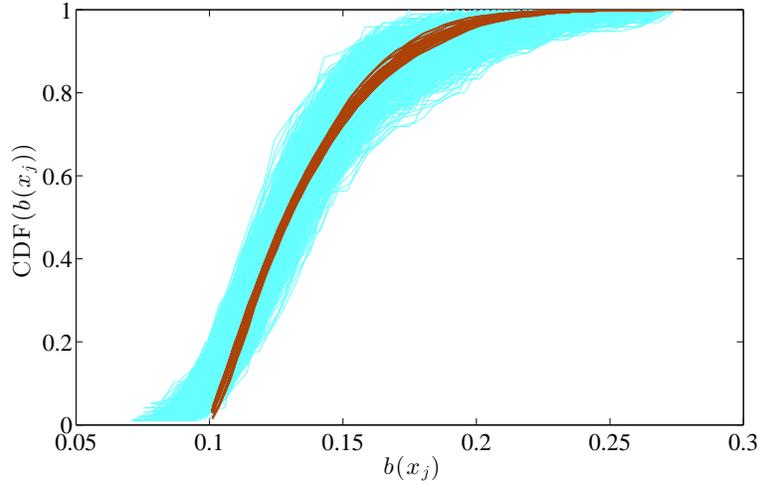}
        \vskip -.3cm
\caption{Thin (turquoise) curves are $N_x$ CDFs of the true field $b(x)$ which form a probability box. Thick (brown) curves are $N_x$ CDFs of the $1$-cut of the empirical field $\tilde{b}(x,{\bf y},{\bf z})$. Each CDF corresponds to a discrete point $\{ x_j \}_{j=1}^{N_x}$ and is obtained from $M=100$ bootstrap samples.}
\label{justify_CDF}
  \end{center}
  \vskip -.5cm
\end{figure}

First, the CDFs of the true field $b(x)$ are obtained as follows. At each observation point $\{ x_j \}_{j=1}^{N_x}$, we use $M$ realizations $\{ b_m(x_j) \}_{m=1}^M$ and find the histogram of $b(x_j)$. We then normalize the histogram and obtain the corresponding CDF. See thin (turquoise) curves in Figure \ref{justify_CDF}, which show $N_x$ CDFs of the true field $b(x)$ corresponding to $N_x$ observation points. Indeed, we obtain a {\it probability box} (or p-box) for the true field $b(x)$, which is another indication of the presence of imprecise uncertainty in the problem. 

Next, we consider the empirical field $\tilde{b}(x,{\bf y},{\bf z})$, constructed by Algorithm \ref{ALG_field}. In the particular case when $L=10^4 [\mu \text{m}]$ and $\ell=100 [\mu \text{m}]$, we need $N=42$ terms in the KL expansion to preserve $85 \%$ of the standard deviation. The random vector ${\bf y}$ therefore contains $N=42$ uncorrelated standard normal random variables. Different $\alpha$-cut intervals of the fuzzy vector ${\bf z}$, obtained in Section \ref{sec:fuzzy_moments}, give different $\alpha$-cuts of the fuzzy-stochastic field $\tilde{b}(x_j,{\bf y},{\bf z})$ at the $N_x$ observation point $\{ x_j \}_{j=1}^{N_x}$. We employ Monte Carlo using $\hat M = 10^4$ independent realizations of the normal random vector ${\bf y}$, and under the assumption that the fuzzy variables are completely interactive, we obtain the CDFs of the left and right limits of $\alpha$-cuts of the empirical field. Thick (brown) curves in Figure \ref{justify_CDF} show $N_x$ CDFs of the $\alpha$-cut for $\alpha=1$. We observe a good agreement between the CDFs of the true field and the CDFs of the constructed empirical field.
%

\subsection{Model Validation}
\label{sec:validation}

The fuzzy-stochastic model constructed in Section \ref{sec:local-model} needs to be validated based on the desired QoI. We note that a fuzzy-stochastic field on $x \in [0,1 \text{m}]$ turns the original problem \eqref{BVP} into a fuzzy-stochastic problem. This will turn the solution \eqref{exact_sol} into a fuzzy-stochastic field, explicitly given by $\tilde{u}(x,{\bf y}, {\bf z}) = \int_{0}^{x} \tilde{b}(\xi,{\bf y}, {\bf z}) \, F(\xi) \, d\xi$. 
%
%

%

In order to validate the proposed model, which is constructed on the small domain $[0,L]$ with $L=10^4 [\mu \text{m}]$, we first need to choose a QoI. For instance, we assume that the same problem  as in \eqref{BVP} holds on the small domain $[0,L]$ with no external force term, i.e. $f \equiv 0$. We consider the solution at the mid-point $x = L/2$ as the QoI, given by:
\begin{equation}\label{QoI1}
{\mathcal Q}_0 = \int_{0}^{L/2} b(\xi) \, d\xi, \qquad 
\tilde{\mathcal Q}_0({\bf y},{\bf z}) = \int_{0}^{L/2} \tilde{b}(\xi,{\bf y},{\bf z}) \, d\xi.
\end{equation}
Here ${\mathcal Q}_0$ is the true QoI computed by the true field $b(x)$, and $\tilde{\mathcal Q}_0$ is the QoI computed by the proposed fuzzy-stochastic model $\tilde{b}(x,{\bf y},{\bf z})$. 
We then proceed with a simple validation strategy \cite{BNT:07_reliability,BNT:08_Sandia}. We compare the true QoI directly obtained by the real data (i.e. $M=100$ bootstrap samples) and the left and right CDFs of the QoI obtained by the proposed model, as follows. First, we find a set of $N_b$ benchmark solutions, referred to as the {\it truth}. For this, we choose $N_b$ groups of samples, where each group consists of $\tilde{M} < M$ different, randomly selected samples out of $M$ bootstrap samples. For each group we then compute $\tilde{M}$ samples of the true quantity ${\mathcal Q}_0$ and then obtain its CDF. This gives us $N_b$ distributions of the true quantity ${\mathcal Q}_0$ in \eqref{QoI1}, which form a p-box and will be used as true validation data. Next, at different $\alpha$-cuts, we compute the CDFs of the left and right limits of the $\alpha$-cuts of the fuzzy-random quantity $\tilde{\mathcal Q}_0$. This can for instance be done by employing a non-intrusive sampling technique, such as Monte Carlo, outlined in Algorithm \ref{ALG_qoi}. 
%
\begin{algorithm}
\caption{Calculating CDFs of left and right limits of $\alpha$-cuts of the fuzzy-random QoI}
\label{ALG_qoi}
\begin{algorithmic} 
\medskip
\STATE {\bf 1.} Select $M_s = 10^4$ independent samples of $N$ unit normal random variables $\{ {\bf y}_m \}_{m=1}^{M_s}$.

\medskip
\STATE {\bf 2.} $\forall \, \alpha$, uniformly discretize the $\alpha$-cut intervals of completely interactive fuzzy variables into $M_{\text{f}}=100$ grid points $\{ {\bf z}_{\alpha k} \}_{k=1}^{M_{\text{f}}}$.  


\medskip
\STATE {\bf 3.} {\bf For} $m=1$ : $M_s$

\STATE {\bf 4.} \hskip .21cm Calculate $\tilde{\mathcal Q}_{0 \alpha l}^m := \min\limits_{k} \tilde{\mathcal Q}_0({\bf y}_m,{\bf z}_{\alpha k})$ and $\tilde{\mathcal Q}_{0 \alpha r}^m := \max\limits_{k} \tilde{\mathcal Q}_0({\bf y}_m,{\bf z}_{\alpha k})$.

\STATE {\bf 5.} {\bf End}

\medskip
\STATE {\bf 6.} From the histograms of $\tilde{\mathcal Q}_{0 \alpha l}$ and $\tilde{\mathcal Q}_{0 \alpha r}$ and find their CDFs.
\end{algorithmic}
\end{algorithm}

\begin{figure}[!ht]
\centering
 \includegraphics[width=0.65\textwidth]{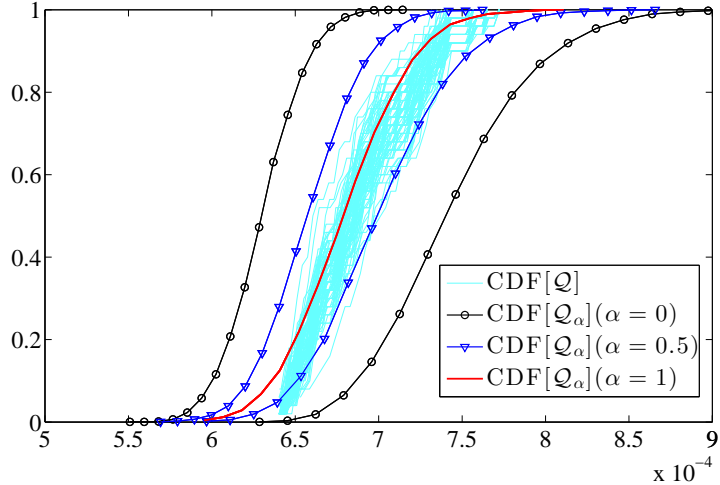}
\vskip -.3cm
\caption{Probability box for the true QoI, directly obtained from the data, and CDFs of the left and right limits of different $\alpha$-cuts of the approximate QoI, computed by the fuzzy-stochastic model.}
\label{probbox_alphacut_H10}
 \end{figure}

In step 4 of Algorith \ref{ALG_qoi}, for a fixed realization of the random vector ${\bf y}_m$, the quantity $\tilde{\mathcal Q}_0$ in \eqref{QoI1} is a fuzzy variable and can be approximated by a quadrature, such as trapezoidal rule: 
\begin{equation} \label{QoI_quad}
\tilde{\mathcal Q}_0({\bf y}_m,{\bf z}) = \int_{0}^{L/2} \tilde{b}(\xi,{\bf y}_m, {\bf z}) \, d\xi \approx \Delta x \,  \sum_{j=1}^{N_M} \tilde{b}(x_j,{\bf y}_m,{\bf z}).
\end{equation}
We need to perform the addition of $N_M$ fuzzy-valued functions $\{ \tilde{b} (x_j,{\bf y}_m,{\bf z}) \}_{j=1}^{N_M}$ of four completely interactive fuzzy variables. The $\alpha$-cuts of the fuzzy quantity in \eqref{QoI_quad} can easily be obtained as described in Section \ref{sec:fuzzy_vec}. 
Figure \ref{probbox_alphacut_H10}, shows the probability box for the true quantity ${\mathcal Q}_0$ and the CDFs of the left and right limits of different $\alpha$-cuts of the approximate quantity $\tilde{\mathcal Q}_0$. 
As we see from the figure, the $\alpha$-cut of the approximate quantity, obtained by the proposed model can accurately capture the p-box of the true quantity. 
It is to be noted that the correlation length $\ell = 100 [\mu \text{m}]$ is actually chosen so that the best fit is obtained.

\section{A global-local multiscale algorithm}
\label{sec:global-local}

In the mechanical design of fiber composites, we often need to determine the local features of the elastic field inside small parts of the domain. For instance, we may need to find the maximum stresses in small zones that are deemed vulnerable to failure. Such problems are amenable to {\it global-local} approaches, in which a {\it global} solution is used to construct a {\it local} solution that captures the microscale features of the true multiscale solution. In this section we develop a global-local technique, based on the proposed fuzzy-stochastic model, to compute the local features of the elastic field. We first start with a short description of the two tools that we need: RVE and homogenization.

\subsection{RVE and homogenization}
\label{sec:RVE-homogenization}

Consider a micro-structure where the size of the domain is sufficiently large relative to the typical inclusion size $h$. For instance, in the present work, the length of the composite plate (1 [m]) is $10^5$ times larger than the typical size of fibers (10 $[\mu \text{m}]$). In general, an RVE is an element which is much smaller than the size of the overall domain and yet contains enough information on the micro-structure (here enough number of fibers) to describe the whole structure, see e. g. \cite{RVE1,RVE2}. Several definitions of an RVE have been proposed based on the above general definition: 
\begin{itemize}
\item RVE is an element that: 1) is structurally entirely typical of the whole mixture on average; and 2) contains a sufficient number of inclusions for the apparent overall moduli to be effectively independent of the surface values of traction and displacement, as long as these values are macroscopically uniform \cite{RVE3}.

\item RVE must be chosen sufficiently large compared to the size of the microstructure for the approach to be valid, and it is the smallest material volume element of the composite for which the usual spatially constant overall macroscopic material properties are a sufficiently accurate model to represent the mean constitutive response \cite{RVE4}.

\item RVE is a model of the material to be used to determine the corresponding effective properties for the homogenized macroscopic model. An RVE should be large enough to contain sufficient information about the microstructure in order to be representative, however it should be much smaller than the macroscopic body. This is known as the Micro-Meso-Macro principle \cite{RVE5}.

\end{itemize}

For micro-structures, an important issue is therefore to find the RVE size. In 1D problems, we need to find the RVE length, denoted by $L_{\text{RVE}}$. By the above definition, the RVE length needs to be obtained from effective (or homogenized) properties. To find the RVE length, we proceed as follows; see also  \cite{RVE1,RVE2}:
\begin{itemize}
\item[1.] Generate a sequence of increasing element lengths $10 \, \mu {\text{m}} \le L_1 < L_2 < \dotsc < L_r \ll 1 \, {\text{m}}$. 
For each $L_i$, with $i=1, \dotsc, r$, we have $M$ samples available, obtained by bootstrapping.  

\item[2.] For each fixed length $L_i$, we use homogenization (see below) with a homogenization length $H=L_i$ and obtain effective parameters for all $M$ samples, denoted by $\{b_{m}^H(x) \}_{m=1}^M$. We then compute the average of sample mean $\bar{\mu}[b^H]$ and the average of sample standard deviation $\bar{\sigma}[b^H]$ of the effective parameter. Note that the sample mean and standard deviation, obtained from $M$ samples, are functions of $x$ given on discrete points along $x$. Their average is taken over all discrete points along $x$.

\item[3.] Set up a criterion and a tolerance, such as:
$$
\varepsilon_{L} := \frac{\bar{\sigma}[b^H]}{\bar{\mu}[b^H]} \le {\text TOL}.
$$
The smallest length $L$ that satisfies the criterion will be chosen as the RVE length $L_{\text{RVE}}$. Note that the larger $L$, or equivalently the larger the homogenization length $H$, the less oscillatory the effective parameters $b_{m}^H(x)$ in $x$, and hence the smaller $\varepsilon_{L}$.
\end{itemize}

\begin{figure}[!h]
\centering
\includegraphics[width=0.7\textwidth]{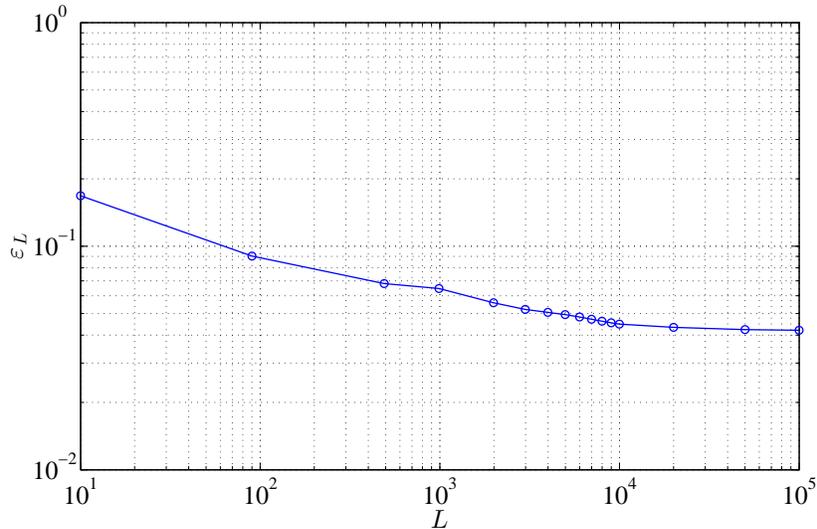}
\caption{An RVE length $L_{\text{RVE}} = 10000 [\mu \text{m}]$ gives a relative error less than $5\%$.}
\label{L_rve}
\end{figure}

Figure \ref{L_rve} shows the relative error $\varepsilon_{L}$ versus $L$. We observe that the length $L = 10^4 [\mu \text{m}]$ gives a relative error less than $5\%$, and the decrease in the error for larger lengths is marginal. We therefore choose $L_{\text{RVE}} = 10^4 [\mu \text{m}]$.

\medskip
\noindent
{\bf Homogenization.} 
A major step in finding the RVE length $L_{\text{RVE}}$ described above is the computation of effective parameters on the RVE. The parameter $a(x)$ is highly oscillatory with variations of the size of fibers ($\sim h = 10 \, \mu$m), which are much smaller than the overall size of the composite ($\sim 1$ m). 
On a given domain of length $L > h$, we may therefore employ periodic homogenization, by periodically extending the highly oscillatory parameters $\{ a_m(x) \}_{m=1}^M$ and obtain (effective) homogenized parameters $\{ a_{m}^H(x) \}_{m=1}^M$, with $H$ being a homogenization length. 

\begin{figure}[!h]
\psfrag{amh10}[r][][0.8]{$a_m^H \, \, \, \, (H=10)$}
\psfrag{amh20}[r][][0.8]{$a_m^H \, \, \, \, (H=20)$}
\psfrag{am}[r][][0.8]{$a_m$}
\psfrag{j1}[c][][0.8]{$j=1$}
\psfrag{j2}[c][][0.8]{$j=2$}
\psfrag{jend}[c][][0.8]{$j=N_x$}
\vskip .3cm
\centering
\includegraphics[width=11cm,height=6cm]{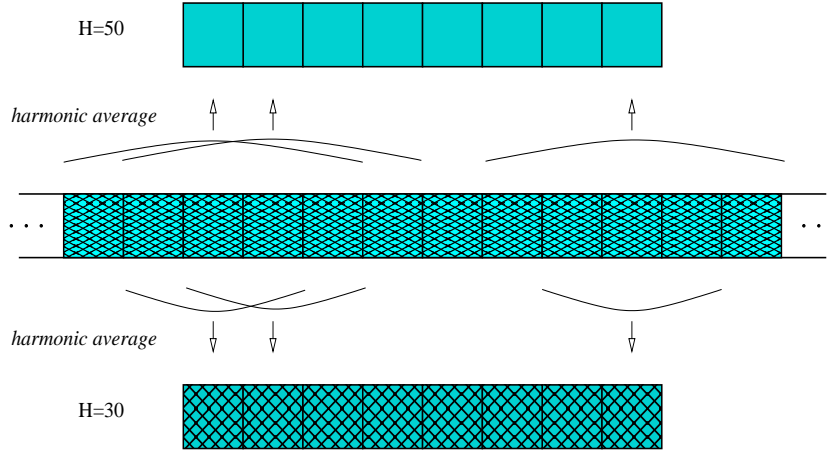}
\caption{A schematic representation of obtaining 1D homogenized data from bootstrap data.}
\label{homogenization_pic}
\end{figure}

Consider a sample $m$ of length $L \ll 1$ m. For any given homogenization length $H$, we find the homogenized parameters $a_{m}^H(x)$ on the sample as follows. We first note that we have the highly oscillatory parameter $a_{m}(x)$ on the whole domain $x \in [0, 1 {\text{m}}]$, given on square elements of size $h \times h$; see Figure \ref{homogenization_pic} (middle).   

The sample with length $L$ is first divided into $N_x$ square elements of size $h \times h$. In each element, we compute the homogenized parameter by taking the harmonic average over the surrounding $H/h$ elements of the sample with highly oscillatory parameters. This gives a homogenized parameter $a_m^H(x)$ on $N_x$ square elements centered at discrete points $\{ x_j \}_{j=1}^{N_x}$. Figure \ref{homogenization_pic} shows the procedure for two different homogenization lengths $H=30, 50$. 
The smallest possible homogenization length in this setting is $H=h$, for which we have $a_m(x) = a_{m}^H(x)$, that corresponds to a full resolution of the micro scale. By repeating the same procedure for all samples $m=1, \dotsc, M$, we collect $M$ realizations of the homogenized parameter $\{a_m^H(x) \}_{m=1}^M$ on a given length $L$. We will eventually use the homogenized parameters $\{b_m^H(x) \}_{m=1}^M$, where $b_m^H = 1/a_m^H$, to find the RVE length $L_{\text{RVE}}$ and to construct the global model.

\subsection{Global-local algorithm}
\label{sec:algorithm}



We will follow a global-local approach and present an algorithm, outlined in Algorithm \ref{ALG_Global_Local}, consisting of two global and local parts. 
\begin{algorithm}
\caption{A global-local algorithm for computing the QoI in \eqref{QOI}} 
\label{ALG_Global_Local}
\begin{algorithmic} 
\medskip

\STATE {\bf 0.} Find $L_{\text{RVE}}$ and select the local domain $D_0 = [x_0 - \frac{L_{\text{RVE}}}{2}, x_0 + \frac{L_{\text{RVE}}}{2}]$ centered at $x_0=0.75$.

\medskip
 
\STATE {\bf Part I.} {\it On global domain $D = [0,1 \text{m}]$:}
\STATE {$\, \,$ \bf 1.} Construct the global (fuzzy-stochastic) field $\tilde{b}_G$ in $D$ (see Section  \ref{sec:global_model}).
\STATE {$\, \,$ \bf 2.} Solve the global problem and find the global solution.

\medskip

\STATE {\bf Part II.} {\it On the local domain $D_0 = [x_0 - \frac{L_{\text{RVE}}}{2}, x_0 + \frac{L_{\text{RVE}}}{2}]$:}
\STATE {$\, \,$ \bf 1.} Compute the boundary data on $\partial D_0$ from the global solution.
\STATE {$\, \,$ \bf 2.} Construct the local (fuzzy-stochastic) field $\tilde{b}_L$ in $D_0$ (see Section  \ref{sec:local_model}).
\STATE {$\, \,$ \bf 3.} Solve the local problem and compute the QoI.
\end{algorithmic}
\end{algorithm}

First, on the global domain $D=[0,1 \text{m}]$, we construct a global model parameter, say $\tilde{b}_{G}$, as described in Section \ref{sec:global_model}, based on the homogenized data on the RVE. We then solve the global problem, which is the same as problem \eqref{BVP} with the $a$ replaced by $\tilde{b}_{G}^{-1}$. The solution to the global problem is used to obtain boundary data for the local domain $D_0 \subset D$, which is the RVE. Inside the local domain $D_0$, we construct a local model parameter, say $\tilde{b}_{L}$, as described in Section \ref{sec:local_model}. Finally, we solve the local problem, which consists of the same equation as \eqref{PDE} with $f \equiv 0$ and Dirichlet boundary conditions with data obtained from the global solution, and compute the QoI.

\subsubsection{Global field}
\label{sec:global_model}

The global field $\tilde{b}_G$ is constructed on the whole domain $D=[0,1 \text{m}]$, and based on the effective parameters on the RVE. The construction is similar to Algorithm \ref{ALG_field} as follows:

\begin{itemize}
\item[1.] Find $L_{\text{RVE}}$ as described in Section \ref{sec:RVE-homogenization}.

\item[2.] Take the homogenization length $H=L_{\text{RVE}}$ and find effective parameters $\{ b_m^H(x_j) \}_{j=1}^{Nx}$ on the RVE for all samples $m=1, \dotsc, M$, as described in Section \ref{sec:RVE-homogenization}. Note that we compute the effective parameters only on the RVE, not on the whole domain $[0, 1 \text{m}]$.

\item[3.] Find the first four moments from the effective samples at all grid points $\{ x_j \}_{j=1}^{N_x}$ and model them by either fuzzy or crisp numbers. Here, since the moments of effective parameters are almost crisp, we model them with crisp values.

\item[4.] Take the correlation length $\ell \sim 5\,  L_{\text{RVE}}$ (see details below).

\item[5.] Construct the empirical global field $\tilde{b}_G$ as in Section \ref{sec:modeling}. The global field will be a non-stationary non-Gaussian translation stochastic field with a small number of KL terms, say $N_{I}$, and slowly varying in $x$. We write ${\tilde b}_G = {\tilde b}_G(x,{\bf y}_I)$, where ${\bf y}_I$ is a vector of $N_I$ uncorrelated standard normal random variables.
\end{itemize}

\begin{figure}[!h]
  \begin{center}
       \subfigure[$H=100$]{\includegraphics[width=8cm,height=4cm]{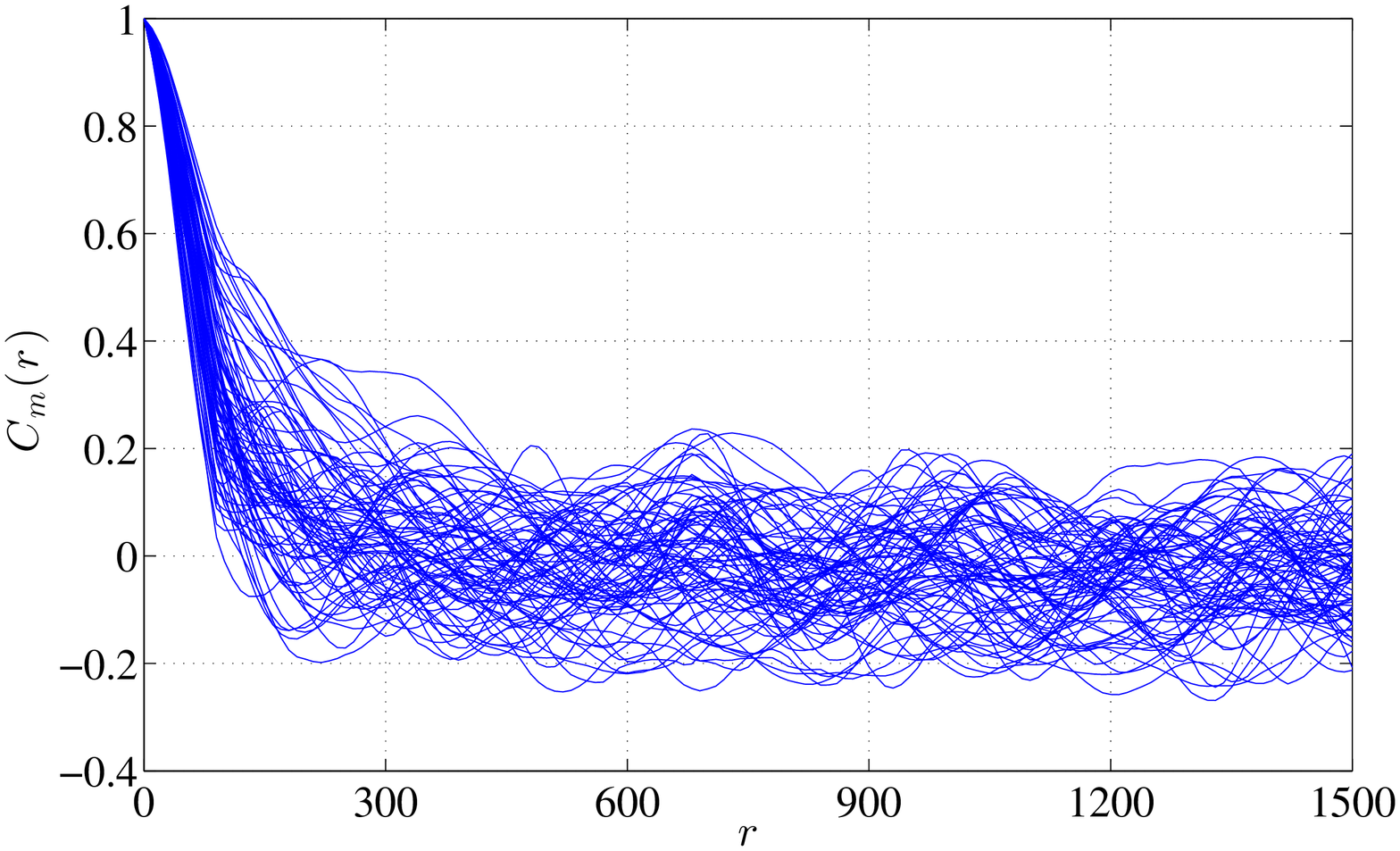}}
        \subfigure[$H=500$]{\includegraphics[width=8cm,height=4cm]{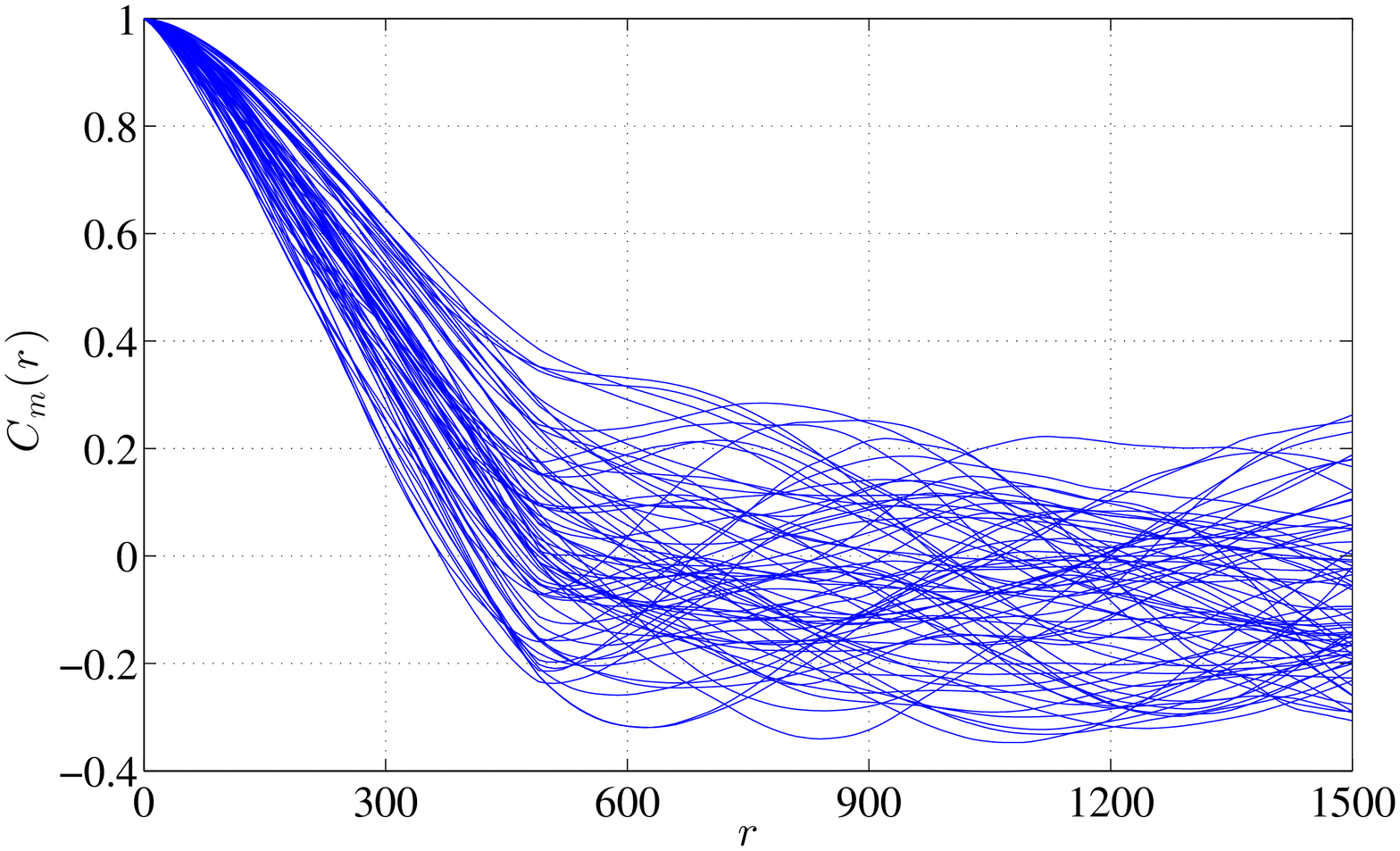}}
\vskip -0.3cm        
\caption{Empirical correlation function $C_m(r_n)$ versus the distance $r_n$ for the homogenized parameter $b_m^H(x)$ with two different $H=100, 500$ and for several realizations $m$.}
\label{CORR2}
  \end{center}
\end{figure}

\noindent
{\bf Correlation length for global field.} 
Figure \ref{CORR2} shows the normalized correlation function $C_m(r_n)$ versus $r_n$, for two different homogenization lengths, $H=100, 500$ and for several realizations $m=1, \dotsc, M$. 
We observe that the correlation function decreases from 1 at $r=0$ by about $50 \%$ at $r=H$, and continues decreasing by $80 \%$ at larger distances. This suggests that the correlation length of the global field $b_G(x)$ is of the order of the homogenization length ($\sim 5 \, H$). We will use this observation to construct the empirical global field. However, we note that the actual value of correlation length for the empirical field will be selected so that the best fit to the data is obtained.

\subsubsection{Local field}
\label{sec:local_model}

The local field $\tilde{b}_L$ is constructed on a local domain of size $L_{\text{RVE}}$ and based on the highly oscillatory (non-homogenized) parameters on the local domain, in the same way as described in Algorithm \ref{ALG_field}. The local field will be a fuzzy-stochastic field $\tilde{b}_L = \tilde{b}_L(x,{\bf y}_{II}, {\bf z})$, where ${\bf y}_{II}$ is a vector of $N_{II}$ uncorrelated standard normal random variables, and ${\bf z}$ is a fuzzy vector of four completely interactive fuzzy variables.

\subsection{Validation}
\label{sec:GL-validation}

Similar to Section \ref{sec:validation}, we need to validate the global-local model and approach. The validation needs to be done based on the true QoI given by \eqref{QOI} and \eqref{exact_sol}:
\begin{equation}\label{true_QOI_GL}
{\mathcal Q} =  \int_{0}^{0.75} b(\xi) \, F(\xi) \, d\xi, \qquad F(\xi) = \left\{ \begin{array}{l l}
2 \, \xi \, & \, \, \xi \in [0 , 0.5),\\
0 \, & \, \, \xi \in [0.5 , 1].
\end{array} \right.
\end{equation}
Here, ${\mathcal Q}_0$ is computed by the true field $b(x)$ on [0, 1 \text{m}]. The approximate QoI, computed by the proposed model and the global-local approach, is obtained following Algorithm \ref{ALG_Global_Local}. The approximate QoI will be a fuzzy-stochastic field, 
\begin{equation}\label{approx_QOI_GL}
\tilde{\mathcal Q} = \tilde{\mathcal Q}({\bf y}, {\bf z}),
\end{equation}
where ${\bf y} = [{\bf y}_I, {\bf y}_{II}]$ is a vector of $N=N_I + N_{II}$ uncorrelated standard normal random variables, and ${\bf z}$ is a vector of four completely interactive fuzzy variables. 

\begin{figure}[!h]
  \begin{center}
      \includegraphics[width=0.65\textwidth]{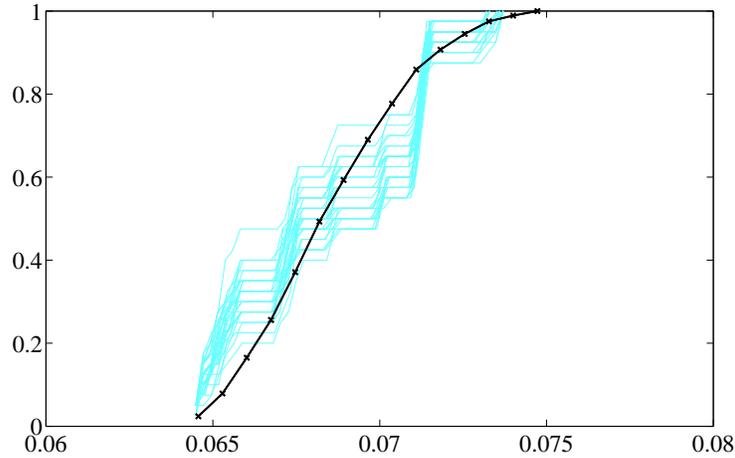}     
\caption{Probability box for the true QoI (thin turquoise curves) and CDF of the 1-cut of the approximate QoI (in black) computed by the fuzzy-stochastic global-local algorithm.}
\label{valid_GL}
  \end{center}
\end{figure}

Similar to the procedure in Section \ref{sec:validation}, we can compare the p-box for the true quantity in \eqref{true_QOI_GL} and the left and right CDFs of the approximate quantity in \eqref{approx_QOI_GL}. As we see in Figure \ref{valid_GL}, the $\alpha$-cut of the approximate quantity, obtained by the proposed method can accurately capture the p-box of the true quantity. 

\section{Conclusions}
\label{sec:CON}

First, we have shown that current multiscale stochastic models, such as stationary random fields, which are based on precise probability theory, are not capable of correctly characterizing uncertainty in fiber composites. Next, we have motivated the applicability of mathematical models based on imprecise uncertainty theory and presented a novel fuzzy-stochastic model, which can more accurately describe uncertainty in fiber composites. The new model combines stochastic fields and fuzzy variables through a simple calibration-validation approach. Finally, we have constructed a global-local multiscale algorithm for efficiently computing output quantities of interest. The algorithm uses the concept of an RVE and computes a global solution to construct a local solution that captures the microscale features of the multiscale solution. The results are based on and backed by real experimental data. 

For simplification and to motivate and establish the main concepts, we have considered a one-dimensional problem in the present work. Future directions include studying fiber composites in two and three dimensions and performing a rigorous error analysis for both the mathematical and computational models. 
%
Application to other engineering problems involving multiple scales and uncertainty is also a subject of our future work.

\section*{Acknowledgements}

The first author is supported by NSF Grant DMS-1211014. The second author would like to recognize the support of the J. T. Oden Fellowship and the Institute for Computational Engineering and Sciences at the University of Texas at Austin.

\bibliographystyle{plain}
\bibliography{refs_me1}

\end{document}